# GLOBAL STABILIZATION OF NONLINEAR DELAY SYSTEMS WITH A COMPACT ABSORBING SET

Iasson Karafyllis[*], Miroslav Krstic[**], Tarek Ahmed-Ali[***],
and Francoise Lamnabhi-Lagarrigue[****]

[*]**Dept. of Mathematics, National Technical University of Athens,**
Zografou Campus, 15780, Athens, Greece (email: iasonkar@central.ntua.gr )

[**]**Dept. of Mechanical and Aerospace Eng., University of California,**
San Diego, La Jolla, CA 92093-0411, U.S.A. (email: krstic@ucsd.edu )

[***]**Laboratoire GREYC CNRS-ENSICAEN**
06 Boulevard du Marechal Juin 14050 Caen Cedex
(email: tarek.ahmed-ali@ensicaen.fr )

[****]**Centre National de la Recherche Scientifique, CNRS-EECI SUPELEC,**
3 rue Joliot Curie, 91192, Gif-sur-Yvette, France
(email: lamnabhi@lss.supelec.fr)

**Abstract**

Predictor-based stabilization results are provided for nonlinear systems with input delays and a compact absorbing set. The control scheme consists of an inter-sample predictor, a global observer, an approximate predictor, and a nominal controller for the delay-free case. The control scheme is applicable even to the case where the measurement is sampled and possibly delayed. The closed-loop system is shown to have the properties of global asymptotic stability and exponential convergence in the disturbance-free case, robustness with respect to perturbations of the sampling schedule, and robustness with respect to measurement errors. In contrast to existing predictor feedback laws, the proposed control scheme utilizes an approximate predictor of a dynamic type which is expressed by a system described by Integral Delay Equations. Additional results are provided for systems that can be transformed to systems with a compact absorbing set by means of a preliminary predictor feedback.

**Keywords:** nonlinear systems, predictor feedback, delay systems.

## 1. Introduction

Remarkable progress has been made in recent years on the design of predictor feedback laws for nonlinear delay systems [3,4,5,7,8,9,10,11,14,15,16,17]. The main challenge to the implementation and design of predictor feedback for nonlinear delay systems is that, except for rare special cases, the solution mapping (used for the prediction) is not available explicitly.

The current status in the literature on input delay compensation is that when, in addition to input delays,

- the full state is not measured,
- the measurement is sampled and possibly delayed,

and when, in addition to global stability, the following properties are required in closed loop,

- exponential convergence for the disturbance-free case,
- robustness with respect to perturbations of the sampling schedule,
- robustness with respect to measurement errors,

predictor feedback designs are available only for two classes of systems: linear detectable and stabilizable systems and globally Lipschitz systems in strict feedback form [10].

In this paper we present a result that removes the global Lipschitz restriction (an algebraic condition on the system's right-hand-side) but imposes an assumption that the system has a compact absorbing set (a condition on the system's dynamic behavior in open loop). Specifically, we consider general nonlinear systems of the form

$$\dot{x}(t) = f(x(t), u(t-\tau)), x \in \Re^n, u \in U \tag{1.1}$$

where $U \subseteq \Re^m$ is a non-empty compact set with $0 \in U$, $\tau \geq 0$ is the input delay and $f : \Re^n \times \Re^m \to \Re^n$ is a smooth vector field with $f(0,0) = 0$. The measurements are sampled and the output is given by

$$y(\tau_i) = h(x(\tau_i - r)) + e(\tau_i) \tag{1.2}$$

where $h : \Re^n \to \Re^k$ is a smooth mapping with $h(0) = 0$, $r \geq 0$ is the measurement delay, $\{\tau_i\}_{i=0}^{\infty}$ is a partition of $\Re_+$ (the set of sampling times) and the input $e : \Re_+ \to \Re$ is the measurement error. We focus on a class of nonlinear systems which is different from the class of globally Lipschitz systems: the *systems with a compact absorbing set*. A nonlinear system with a compact absorbing set is a system for which all solutions enter a specific compact set after an initial transient period (for systems without inputs the name "global uniform ultimate boundedness" is used in [13]; the term "dissipative system" is used in the literature of finite-dimensional dynamical systems; see [20] and the discussion on page 22 of the book [21]).

Though it may appear that we merely trade one major restriction (global Lipschitzness) for another (compact absorbing set), which imposes a strong requirement on the system's open-loop behavior, the latter restriction is less frequently violated in applications. Many engineering systems belong to the class of systems with a compact absorbing set because finite escape is rare in physical processes, control inputs usually saturate, and limit cycles are a frequent outcome of local instabilities.

The contribution of our paper is twofold:

a) predictor feedback is designed and stability is proved for the class of nonlinear delay systems with a compact absorbing set under appropriate assumptions (Theorem 2.2),
b) the result is then extended to nonlinear delay systems that can be transformed to systems with a compact absorbing set by means of a preliminary predictor feedback (Theorem 2.4).

In both cases, we provide explicit formulae for the predictor feedback and explicit inequalities for the parameters of the applied control scheme and the upper diameter of the sampling partition. The proposed predictor feedback guarantees all properties listed at the beginning of the section for the class of nonlinear delay systems with a compact absorbing set: global asymptotic stability and global exponential attractivity in the absence of measurement error, robustness with respect to perturbations of the sampling schedule and robustness with respect to measurement errors.

Our predictor feedback design consists of the following elements:
1) an Inter-Sample Predictor (ISP), which uses the sampled, delayed and corrupted measurements of the output and provides an estimate of the (unavailable) delayed continuous output signal,
2) a global observer (O), which uses the estimate of the delayed continuous output signal and provides an estimate of the delayed state vector,
3) an approximate or exact predictor (P), which uses the estimate of the delayed state vector in order to provide an estimate of the future state vector, and
4) a delay-free controller (DFC), i.e., a baseline feedback law that works for the delay-free version of the system, which in the presence of delay uses the estimate of the future state vector in order to provide the control action.

We refer to the above control scheme as the ISP-O-P-DFC control scheme. In [10] the ISP-O-P-DFC control scheme was shown to achieve all the objectives mentioned at the beginning of this section by using approximate predictors that are based on *successive approximations of the solution map* for linear detectable and stabilizable systems and globally Lipschitz systems in strict feedback form. Here, we show that the ISP-O-P-DFC control scheme guarantees all the objectives listed at the beginning of this section using *dynamic approximate predictors* for systems with a compact absorbing set.

This methodological difference relative to [10] merits further emphasis. We employ here a class of approximate predictors that are implemented by means of a dynamical system: the approximate predictor is a system described by Integral Delay Equations (IDEs; see [12]) and consists of a series connection of $N$ approximate predictors (each making a prediction for the state vector $\delta = \frac{r+\tau}{N}$ time units ahead). Such dynamic predictors were introduced in [2,6] but here the predictor is designed in a novel way so that the prediction takes values in an appropriate compact set after an initial transient period. The dynamic predictor is different from other predictors proposed in the literature (e.g., exact predictors in [9,17], approximate predictors based on successive approximations in [7,10], approximate predictors based on numerical schemes in [11]). Theorem 2.4 employs a novel combination of approximate predictors and exact predictors in the control scheme, which can be used for other classes of nonlinear delay systems.

The main advantage of the dynamic predictor employed here over other predictor approximations (numerical [11] or successive approximations [7,10]) is the existence of simple formulas (provided in [2]), for the estimation of the asymptotic gain of the measurement error for certain classes of systems. In contrast, the predictor for which the effect of measurement errors is most difficult to quantify is the numerical predictor [11].

On the other hand, the disadvantages of the dynamic predictor are the difficulty of implementation (one has to approximate numerically the solution of the IDEs or the equivalent distributed delay differential equations) and that it works only for certain classes of nonlinear systems (globally Lipschitz systems and systems with a compact absorbing set). In contrast, the most easily programmable predictor is the numerical predictor [11], which is the crudest version of the predictor based on successive approximations [7,10] -- when only one successive approximation is used (and many grid points), then the predictor based on successive approximations coincides with the numerical predictor.

Though our approach to stabilization of nonlinear systems with actuation and measurement delays is based on delay compensation via predictor design—an approach known for its ability to recover nominal performance in the absence of delay and after finite time in the presence of delay—this is not the only option for stabilization of nonlinear systems with large dead times. For certain classes of nonlinear systems other approaches exist that are capable of guaranteeing stability and robustness [18,19].

The paper is structured as follows. Section 2 contains the assumptions and the statements of the main results. The proofs are given in Section 3. Section 4 presents two illustrative examples. Concluding remarks are provided in Section 5.

**Notation.** Throughout this paper, we adopt the following notation:

* $\Re_+ := [0,+\infty)$. A partition of $\Re_+$ is an increasing sequence $\{\tau_i\}_{i=0}^\infty$ with $\tau_0 = 0$ and $\lim_{i \to \infty} \tau_i = +\infty$.
* By $C^0(A;\Omega)$, we denote the class of continuous functions on $A \subseteq \Re^n$, which take values in $\Omega \subseteq \Re^m$. By $C^k(A;\Omega)$, where $k \geq 1$ is an integer, we denote the class of functions on $A \subseteq \Re^n$ with continuous derivatives of order $k$, which take values in $\Omega \subseteq \Re^m$.
* By $\text{int}(A)$, we denote the interior of the set $A \subseteq \Re^n$.
* For a vector $x \in \Re^n$, we denote by $x'$ its transpose and by $|x|$ its Euclidean norm. $A' \in \Re^{n \times m}$ denotes the transpose of the matrix $A \in \Re^{m \times n}$ and $|A|$ denotes the induced norm of the matrix $A \in \Re^{m \times n}$, i.e., $|A| = \sup\{|Ax| : x \in \Re^m, |x| = 1\}$.
* A function $V : \Re^n \to \Re_+$ is called positive definite if $V(0) = 0$ and $V(x) > 0$ for all $x \neq 0$. A function $V : \Re^n \to \Re_+$ is called radially unbounded if the sets $\{x \in \Re^n : V(x) \leq M\}$ are either empty or bounded for all $M \geq 0$.
* For a function $V \in C^1(A;\Re)$, the gradient of $V$ at $x \in A \subseteq \Re^n$, denoted by $\nabla V(x)$, is the row vector $\nabla V(x) = \left[ \frac{\partial V}{\partial x_1}(x) \quad \ldots \quad \frac{\partial V}{\partial x_n}(x) \right]$.
* The class of functions $K_\infty$ is the class of strictly increasing, continuous functions $a : \Re_+ \to \Re_+$ with $a(0) = 0$ and $\lim_{s \to +\infty} a(s) = +\infty$.

## 2. Systems with an Absorbing Compact Set

Consider the system (1.1), (1.2). Our main assumption guarantees that there exists a compact set which is robustly globally asymptotically stable (the adjective *robust* means uniformity to all measurable and essentially bounded inputs $u : \Re_+ \to U$). We call the compact set "absorbing" because the solution "is absorbed" in the set after an initial transient period.

**(H1)** *There exist a radially unbounded (but not necessarily positive definite) function $V \in C^2(\Re^n; \Re_+)$, a positive definite function $W \in C^1(\Re^n; \Re_+)$ and a constant $R > 0$ such that the following inequality holds for all $(x,u) \in \Re^n \times U$ with $V(x) \geq R$*

$$\nabla V(x) f(x,u) \leq -W(x) \qquad (2.1)$$

Indeed, assumption (H1) guarantees that for every initial condition $x(0) \in \Re^n$ and for every measurable and essentially bounded input $u : \Re_+ \to U$ the solution $x(t)$ of (1.1) enters the compact set $S = \{x \in \Re^n : V(x) \leq R\}$ after a finite transient period, i.e., there exists $T \in C^0(\Re^n; \Re_+)$ such that $x(t) \in S$, for all $t \geq T(x(0))$. Moreover, notice that the compact set $S = \{x \in \Re^n : V(x) \leq R\}$ is positively invariant. This fact is guaranteed by the following lemma which is an extension of Theorem 5.1 in [13] (page 211).

**Lemma 2.1:** *Consider system (1.1) under hypothesis (H1). Then there exists $T \in C^0(\Re^n; \Re_+)$ such that for every $x_0 \in \Re^n$ and for every measurable and essentially bounded input $u : [-\tau, +\infty) \to U$ the solution $x(t) \in \Re^n$ of (1.1) with initial condition $x(0) = x_0$ and corresponding to input $u : [-\tau, +\infty) \to U$ satisfies $V(x(t)) \leq \max(V(x_0), R)$ for all $t \geq 0$ and $V(x(t)) \leq R$ for all $t \geq T(x_0)$.*

Our second assumption guarantees that we are in a position to construct an appropriate local exponential stabilizer for the delay-free version system (1.1), i.e., system (1.1) with $\tau = 0$.

**(H2)** *There exist a positive definite function $P \in C^2(\Re^n; \Re_+)$, constants $\mu, K_1 > 0$ with $K_1|x|^2 \leq P(x)$ for all $x \in \Re^n$ with $V(x) \leq R$ and a globally Lipschitz mapping $k : \Re^n \to U$ with $k(0) = 0$ such that the following inequality holds*

$$\nabla P(x) f(x, k(x)) \leq -2\mu |x|^2, \text{ for all } x \in \Re^n \text{ with } V(x) \leq R \tag{2.2}$$

The requirement that the mapping $k : \Re^n \to U$ with $k(0) = 0$ is globally Lipschitz is not essential. Notice that if assumption (H2) holds for certain locally Lipschitz mapping $k : \Re^n \to U$ with $k(0) = 0$ then we are in a position to define $\tilde{k}(x) = (R + 1 - \min(R + 1, \max(R, V(x))))k(x)$ and we notice that assumption (H2) holds for the globally Lipschitz function $\tilde{k} : \Re^n \to U$.

Our third assumption guarantees that we are in a position to construct an appropriate local exponential observer for the delay-free version of system (1.1), (1.2), i.e., system (1.1), (1.2) with $r = \tau = 0$.

**(H3)** *There exist a symmetric and positive definite matrix $Q \in \Re^{n \times n}$, constants $\omega > 0$, $b > R$ and a matrix $L \in \Re^{n \times k}$ such that the following inequality holds*

$$(z-x)'Q(f(z,u) + L(h(z) - h(x)) - f(x,u)) \leq -\omega |z-x|^2,$$
$$\text{for all } u \in U, z, x \in \Re^n \text{ with } V(z) \leq b \text{ and } V(x) \leq R \tag{2.3}$$

Indeed, assumption (H3) in conjunction with assumption (H1) guarantees that for every $x(0) \in S = \{x \in \Re^n : V(x) \leq R\}$ and for every measurable and essentially bounded input $u : \Re_+ \to U$ the solution of system (1.1), (1.2) with

$$\dot{z} = f(z, u) + L(h(z) - y) \tag{2.4}$$

satisfies an estimate of the form $|z(t) - x(t)| \leq M \exp(-\sigma t)|z(0) - x(0)|$, for all $t \geq 0$ for appropriate constants $M, \sigma > 0$, provided that the initial estimation error $|z(0) - x(0)|$ is sufficiently small. This is why system (2.4) is termed "a local exponential observer". The reader should notice that assumption (H3) holds automatically for nonlinear systems of the form

$$\begin{aligned}
\dot{x}_1 &= f_1(x_1) + x_2 \\
\dot{x}_2 &= f_2(x_1, x_2) + x_3 \\
&\vdots \\
\dot{x}_n &= f_n(x_1,\ldots,x_n) + u \\
y &= x_1
\end{aligned} \qquad (2.5)$$

for every $b > R > 0$ and for every non-empty set $U \subseteq \Re^m$, where $f_i : \Re^i \to \Re$ ($i = 1,\ldots,n$) are smooth mappings.

In order to be able to construct a feedback stabilizer for system (1.1), (1.2) we need an additional technical assumption.

**(H4)** *There exist constants $c \in (0,1)$, $R \leq a < b$ such that the following inequality holds:*

$$\nabla V(z)(f(z,u) + L(h(z) - h(x))) \leq -W(z) + (1-c)|\nabla V(z)|^2 \frac{(z-x)'Q(f(z,u) + L(h(z) - h(x)) - f(x,u))}{\nabla V(z)Q(z-x)}$$

*for all $u \in U$, $z, x \in \Re^n$ with $a < V(z) \leq b$, $\nabla V(z)Q(z-x) < 0$ and $V(x) \leq R$* (2.6)

Assumption (H4) imposes constraints for the evolution of the trajectories of the local observer (2.4). Indeed, inequality (2.6) imposes a bound on the derivative of the Lyapunov function $V \in C^1(\Re^n; \Re_+)$ along the trajectories of the local observer (2.4) for certain regions of the state space.

We are now ready to state the first main result of the paper. Notice that the dynamic feedback stabilizer is explicitly given and that all parameters included in the feedback stabilizer are requited to satisfy explicit inequalities that can be verified easily in practice.

**Theorem 2.2:** *Consider system (1.1), (1.2) under assumptions (H1-4). Define:*

$$\hat{k}(z,y,u) := L(h(z) - y), \text{ for all } (z,y,u) \in \Re^n \times \Re^k \times U \text{ with } V(z) \leq R \qquad (2.7)$$

$$\hat{k}(z,y,u) := L(h(z) - y) - \frac{\varphi(z,y,u)}{|\nabla V(z)|^2}(\nabla V(z))', \text{ for all } (z,y,u) \in \Re^n \times \Re^k \times U \text{ with } V(z) > R \qquad (2.8)$$

*where $\varphi : \Re^n \times \Re^k \times \Re^m \to \Re_+$ is defined by*

$$\varphi(z,y,u) := \max\bigl(0, \nabla V(z)f(z,u) + W(z) + p(V(z))\nabla V(z)L(h(z) - y)\bigr) \qquad (2.9)$$

*and $p : \Re_+ \to [0,1]$ is an arbitrary locally Lipschitz function that satisfies $p(s) = 1$ for all $s \geq b$ and $p(s) = 0$ for all $s \leq a$. Let $q : \Re \to \Re_+$ be a continuously differentiable function with $q(s) = 1$ for $s \leq 1$ and $sq(s) \leq K$ for $s \geq 1$, where $K \geq 1$ is a constant. Let $\psi : \Re^n \to [1,+\infty)$ be a smooth function that satisfies the following implication:*

$$V(x) \leq \max(V(z), b) \Rightarrow |x| \leq \psi(z) \qquad (2.10)$$

*Let $N > 0$ be an integer and $T_s > 0$, $\sigma > 0$ be constants so that:*

$$\sigma \leq \min\left(\frac{\mu}{\sqrt{n}\widetilde{P}}, \frac{c\omega}{4|Q|}\right), \quad \delta M_1^q M_1^f e^{\sigma\delta} < 1 \quad \text{and} \quad T_s G_1 e^{\sigma T_s} \sqrt{\frac{2|Q|}{K_2}} \frac{G_2|Q|}{c\omega} < 1 \tag{2.11}$$

where $M_1^q := \sup\left\{\frac{\left|q\left(\frac{|\xi|}{\psi(z)}\right)\xi - q\left(\frac{|x|}{\psi(z)}\right)x\right|}{|x-\xi|} : x \in S_1, \xi \in S_3, z \in S_2, x \neq \xi\right\}$, $M_1^f := \sup\left\{\frac{|f(x,u)-f(\xi,u)|}{|x-\xi|} : x \in S_1, \xi \in S_4, u \in U, x \neq \xi\right\}$,

$S_1 := \{x \in \Re^n : V(x) \leq R\}$, $S_2 := \{x \in \Re^n : V(x) \leq b\}$, $S_3 := \{x \in \Re^n : |x| \leq K\psi(z) + \delta p(K\psi(z)), z \in S_2\}$, $\delta = \frac{r+\tau}{N}$,

$G_1 := \sup\left\{\frac{|\nabla h(x)f(x,u) - \nabla h(z)f(z,u)|}{|x-z|} : x \in S_1, z \in S_2, u \in U, x \neq z\right\}$, $p(s) := \max\{|f(z,u)| : (z,u) \in \Re^n \times U, |z| \leq s\}$,

$G_2 := \sup\left\{\frac{|\hat{k}(z,y,u) - \hat{k}(z,w,u)|}{|y-w|} : y, w \in \Re^k, z \in S_2, u \in U, y \neq w\right\}$, $S_4 := \{x \in \Re^n : |x| \leq K\psi(z), z \in S_2\}$,

$\widetilde{P} := \max\{|\nabla^2 P(x)| : x \in co(S_1)\}$, $co(S_1)$ *denotes the convex hull of* $S_1$ *and* $K_2 \in (0, |Q|]$ *is a constant for which the inequality* $K_2|x|^2 \leq x'Qx$ *for all* $x \in \Re^n$.

*Then there exist a constant* $\Gamma > 0$ *and a locally Lipschitz function* $C \in K_\infty$ *such that for every partition* $\{\tau_i\}_{i=0}^\infty$ *of* $\Re_+$ *with* $\sup_{i \geq 0}(\tau_{i+1} - \tau_i) \leq T_s$, $e \in L^\infty(\Re_+; \Re^k)$, $\xi_{i,0} \in L^\infty([-\delta, 0); \Re^n)$ ($i = 1, \ldots, N$), $(z_0, w_0) \in \Re^n \times \Re^k$, $x_0 \in C^0([-r, 0]; \Re^n)$, $u_0 \in L^\infty([-r-\tau, 0); U)$, *the solution of (1.1), (1.2) with*

$$\dot{z}(t) = f(z(t), u(t-r-\tau)) + \hat{k}(z(t), w(t), u(t-r-\tau)), \text{ for } t \geq 0 \tag{2.12}$$

$$\dot{w}(t) = \nabla h(z(t))f(z(t), u(t-r-\tau)), \text{ for } t \in [\tau_i, \tau_{i+1}), i \geq 0 \tag{2.13}$$

$$w(\tau_i) = y(\tau_i), \text{ for } i \geq 1 \tag{2.14}$$

$$\xi_j(t) = q\left(\frac{|\xi_{j-1}(t)|}{\psi(z(t))}\right)\xi_{j-1}(t) + \int_0^\delta f\left(q\left(\frac{|\xi_j(t+s-\delta)|}{\psi(z(t))}\right)\xi_j(t+s-\delta), u(t+(j-1)\delta - \tau - r + s)\right)ds, \text{ for } t \geq 0, j = 1, \ldots, N \tag{2.15}$$

*with* $\xi_0(t) = z(t)$ *and*

$$u(t) = k(\xi_N(t)), \text{ for } t \geq 0 \tag{2.16}$$

*initial condition* $\xi_j(\theta) = \xi_{j,0}(\theta)$ *for* $\theta \in [-\delta, 0)$ ($j = 1, \ldots, N$), $(z(0), w(0)) = (z_0, w_0)$, $x(\theta) = x_0(\theta)$ *for* $\theta \in [-r, 0]$, $u(\theta) = u_0(\theta)$ *for* $\theta \in [-r-\tau, 0)$, *exists and satisfies the following estimate for all* $t \geq 0$:

$$\sup_{-r \leq s \leq 0}(|x(t+s)|) + |w(t)| + |z(t)| + \sum_{j=1}^N \sup_{-\delta \leq s < 0}(|\xi_j(t+s)|) + \sup_{-r-\tau \leq s < 0}(|u(t+s)|) \leq$$
$$e^{-\sigma t}C\left(\sup_{-r \leq s \leq 0}(|x_0(s)|) + |z_0| + |w_0| + \sum_{j=1}^N \sup_{-\delta \leq s < 0}(|\xi_{j,0}(s)|) + \sup_{-r-\tau \leq s < 0}(|u_0(s)|) + \sup_{0 \leq s \leq t}(|e(s)|)\right) + \Gamma \sup_{0 \leq s \leq t}(|e(s)|) \tag{2.17}$$

**Remark 2.3:**

**(a)** As noted in the Introduction, Theorem 2.2 shows that the control scheme which consists of the series connection of (i) the sampled-data hybrid observer (2.12), (2.13), (2.14) which provides an estimate of the delayed state vector $x(t-r)$, (ii) the dynamic predictor (2.15) which provides an estimate of the future value of the state vector $x(t+\tau)$ and (iii) the control law (2.16) is successful under assumptions (H1-4) provided that the upper diameter of the sampling partition $T_s > 0$ is sufficiently small. The result of Theorem 2.2 guarantees robustness with respect to perturbations of the sampling schedule (inequality (2.17) holds for every sampling partition $\{\tau_i\}_{i=0}^{\infty}$ with upper diameter less or equal to $T_s > 0$).

**(b)** The sampled-data hybrid observer (2.12), (2.13), (2.14) uses the local exponential observer involved in assumptions (H3), (H4) with some modifications. The first modification involves the replacement of the unavailable output signal $y(t-r)$ with the signal $w(t)$, which is generated by the intersample predictor (2.13), (2.14) (see also [2]). The second modification is the addition of a "correction term" of the form $-\frac{\varphi(z,y,u)}{|\nabla V(z)|^2}(\nabla V(z))'$ which has the task to guarantee the validity of the differential inequality $\nabla V(z)(f(z,u) + \hat{k}(z,y,u)) \leq -W(z)$ for all $(z,y,u) \in \Re^n \times \Re^k \times U$ with $V(z) \geq b$. The "correction term" $-\frac{\varphi(z,y,u)}{|\nabla V(z)|^2}(\nabla V(z))'$ was used in [2] in order to guarantee that the solution enters an appropriate compact set in finite time and in this appropriate compact set the local exponential observer works.

**(c)** The input $e: \Re_+ \to \Re$ quantifies the effect of measurement errors. Inequality (2.17) shows that the "asymptotic gain" of the closed-loop system with respect to the measurement error is linear, i.e. $\limsup_{t \to +\infty} \left( \sup_{-r \leq s \leq 0}(|x(t+s)|) + |w(t)| + |z(t)| + \sum_{j=1}^{N} \sup_{-\delta \leq s < 0}(|\xi_j(t+s)|) + \sup_{-r-\tau \leq s < 0}(|u(t+s)|) \right) \leq \Gamma \limsup_{t \to +\infty} |e(t)|$.
However, the ISS-like inequality (2.17) does not guarantee the ISS property with linear gain. In general, the locally Lipschitz function $C \in K_\infty$ is nonlinear and the gain function with respect to the measurement error is nonlinear.

**(d)** The predictor (2.15) is a system described by Integral Delay Equations (IDEs; see [12]) and consists of the series connection of $N$ predictors (each making a prediction for the state vector $\delta$ time units ahead). Such dynamic predictors were used in [2,6] but here the predictor (2.15) has an important difference with other predictors: the use of the terms $q\left(\frac{|\xi_{j-1}(t)|}{\psi(z(t))}\right)\xi_{j-1}(t)$ instead of $\xi_{j-1}(t)$ guarantees that the prediction will take values in an appropriate compact set. The dynamic predictor (2.15) is different from other predictors proposed in the literature (e.g., exact predictors in [9,17], approximate predictors based on successive approximations in [7,10], approximate predictors based on numerical schemes in [11]).

**(e)** An example of a function $q: \Re \to \Re_+$ that satisfies the requirements of Theorem 2.2 is the function $q(s) := 2s^{-1} - s^{-2}$ for $s > 1$ and $q(s) = 1$ for $s \leq 1$.

**(f)** Since the function $C \in K_\infty$ is a locally Lipschitz function, it follows from estimate (2.17) that the dynamic hybrid controller (2.12), (2.13), (2.14), (2.15) and (2.16) guarantees not only global asymptotic stability but local exponential stability as well in the absence of measurement error. Notice that the stability properties of the closed-loop system are robust with respect to perturbations of the sampling schedule.

The following result uses a preliminary predictor feedback in order to transform the given system to a system with a compact absorbing set. However, the result of Theorem 2.4 does not allow us to conclude Robust Global Asymptotic Stability for the closed-loop system: only exponential attractivity holds for the closed-loop system. The notion of forward completeness used in the statement of Theorem 2.4 is the standard notion used in [1]: the solution exists for times, all initial conditions and all measurable and locally essentially bounded inputs.

**Theorem 2.4:** *Consider the forward complete system*

$$\dot{x}(t) = \tilde{f}(x(t), v(t-\tau))$$
$$x(t) \in \Re^n, v(t) \in \Re^m$$
(2.18)

*where $\tau > 0$ is the input delay, $\tilde{f}: \Re^n \times \Re^m \to \Re^n$ is a smooth vector field with $\tilde{f}(0,0) = 0$ and sampled measurements given by*

$$y(\tau_i) = h(x(\tau_i - r))$$
(2.19)

*where $h: \Re^n \to \Re^k$ is a smooth mapping with $h(0) = 0$, $r \geq 0$ is the measurement delay and $\{\tau_i\}_{i=0}^{\infty}$ is a partition of $\Re_+$ (the set of sampling times).*

*Suppose that there exist smooth functions $a_1: \Re^n \to \Re^l$, $a_2: \Re^l \to \Re^m$ with $a_1(0) = 0$, $a_1(0) = 0$ and a non-empty compact set $U \subset \Re^m$ with $0 \in U$ such that the vector field $f(x,u) := \tilde{f}(x, a_2(a_1(x)) + u)$ satisfies assumptions (H1-4). Moreover, suppose that there exists a locally Lipschitz vector field $g: \Re^l \times \Re^m \to \Re^l$ with $g(0,0) = 0$ such that the equation $\nabla a_1(x)\tilde{f}(x,v) = g(a_1(x), v)$ holds for all $(x,v) \in \Re^n \times \Re^m$. Assume that the system*

$$\dot{\theta}(t) = g(\theta(t), a_2(\theta(t)) + u(t))$$
$$\theta(t) \in \Re^l, u(t) \in U$$
(2.20)

*is forward complete. Finally, suppose that there exists a mapping $\Phi: \Re^k \times L^\infty([-\tau-r,0]; \Re^m) \to$ such that for every $(x_0, v) \in \Re^n \times L^\infty([-\tau, +\infty); \Re^m)$ the solution $x(t) \in \Re^n$ of (2.18) with initial condition $x(0) = x_0$ corresponding to input $v \in L^\infty([-\tau, +\infty); \Re^m)$ satisfies for all $t \geq r$:*

$$a_1(x(t+\tau)) = \Phi(h(x(t-r)), v_t)$$
(2.21)

*where $(v_t)(s) = v(t+s)$ for $s \in [-\tau-r, 0]$.*

*Let $\hat{k}: \Re^n \times \Re^k \times \Re^m \to \Re^n$ be the vector field defined by (2.7), (2.8), (2.9) for certain locally Lipschitz function $p: \Re_+ \to [0,1]$ that satisfies $p(s) = 1$ for all $s \geq b$ and $p(s) = 0$ for all $s \leq a$. Let $q: \Re \to \Re_+$ be a continuously differentiable function with $q(s) = 1$ for $s \leq 1$ and $sq(s) \leq K$ for $s \geq 1$, where $K \geq 1$ is a constant. Let $\psi: \Re^n \to [1, +\infty)$ be a smooth function that satisfies implication (2.10). Let $N > 0$ be an integer and $T_s > 0$, $\sigma > 0$ be constants so that (2.11) holds.*

*Then for every partition $\{\tau_i\}_{i=0}^{\infty}$ of $\Re_+$ with $\sup_{i \geq 0}(\tau_{i+1} - \tau_i) \leq T_s$, $\xi_{i,0} \in L^\infty([-\delta, 0]; \Re^n)$ ($i = 1, ..., N$), $(z_0, w_0, \theta_0) \in \Re^n \times \Re^k \times \Re^l$, $x_0 \in C^0([-r, 0]; \Re^n)$, $u_0 \in L^\infty([-r-\tau, 0]; U)$, $v_0 \in L^\infty([-\tau, 0]; \Re^m)$, the solution of (2.18), (2.19) with (2.12), (2.13), (2.14), (2.15), (2.16) and*

$$v(t) = u(t) + a_2(\theta(t)), \text{ for } t \geq 0 \tag{2.22}$$

$$\dot{\theta}(t) = g(\theta(t), a_2(\theta(t)) + u(t)), \text{ for } t \in [\tau_i, \tau_{i+1}), \ i \geq 0 \tag{2.23}$$

$$\theta(\tau_i) = \Phi(y(\tau_i), v_{\tau_i}), \text{ for } i \geq 1 \tag{2.24}$$

*with* $\xi_0(t) = z(t)$ *and initial condition* $\xi_j(s) = \xi_{j,0}(s)$ *for* $s \in [-\delta, 0)$ ($j = 1, ..., N$), $(z(0), w(0), \theta(0)) = (z_0, w_0, \theta_0)$, $x(s) = x_0(s)$ *for* $s \in [-r, 0]$, $u(s) = u_0(s)$ *for* $s \in [-r-\tau, 0)$, $v(s) = v_0(s)$ *for* $s \in [-\tau, 0)$ *exists for all* $t \geq 0$ *and satisfies:*

$$\limsup_{t \to +\infty} \left( e^{\sigma t} P(t) \right) < +\infty \tag{2.25}$$

*where* $P(t) := \sup_{-r \leq s \leq 0} \left( |x(t+s)| \right) + |w(t)| + |z(t)| + |\theta(t)| + \sum_{j=1}^{N} \sup_{-\delta \leq s < 0} \left( |\xi_j(t+s)| \right) + \sup_{-r-\tau \leq s < 0} \left( |u(t+s)| \right).$

**Remark 2.5:** Theorem 2.4 uses a combination of exact predictors (in the spirit of [9]; (2.21) is an exact prediction of $a_1(x(t+\tau))$ ) and approximate predictors ((2.15), (2.16) provides an approximate prediction of $k(x(t+\tau))$ ). Therefore, Theorem 2.4 generalizes the results in [9] and the result of Theorem 2.2. However, as remarked above the result of Theorem 2.4 is simple exponential attractivity for the closed-loop system. The existence of functions $a_1 : \Re^n \to \Re^l$, $a_2 : \Re^l \to \Re^m$ satisfying the assumptions of Theorem 2.4 is a restrictive assumption, which can be verified in certain cases (see Example 4.2 below).

## 3. Proofs of Main Results

We start with the proof of Theorem 2.2.

**Proof of Theorem 2.2:** We first notice that the following inequality holds for all $(z, w, u) \in \Re^n \times \Re^k \times U$ with $V(z) \geq b$:

$$\nabla V(z)(f(z,u) + \hat{k}(z,w,u)) \leq -W(z) \tag{3.1}$$

Definition (2.8) implies $\nabla V(z)(f(z,u) + \hat{k}(z,w,u)) = \nabla V(z)(f(z,u) + L(h(z)-w)) - \varphi(z,w,u)$. By distinguishing the cases $\nabla V(z) f(z,u) + W(z) + \nabla V(z) L(h(z)-w) \leq 0$ and $\nabla V(z) f(z,u) + W(z) + \nabla V(z) L(h(z)-w) > 0$, using definition (2.9) and noticing that $p(V(z)) = 1$ we conclude that (3.1) holds.

Let $\{\tau_i\}_{i=0}^{\infty}$ be a partition of $\Re_+$ with $\sup_{i \geq 0}(\tau_{i+1} - \tau_i) \leq T_s$, $x_0 \in C^0\left([-r,0]; \Re^n\right)$, $u_0 \in L^{\infty}\left([-r-\tau,0); U\right)$, $\xi_{i,0} \in L^{\infty}\left([-\delta, 0); \Re^n\right)$ ($i = 1, ..., N$), $(z_0, w_0) \in \Re^n \times \Re^k$, $e \in L^{\infty}_{loc}\left(\Re_+; \Re^k\right)$ and consider the solution of (1.1), (2.12), (2.13), (2.14), (2.15), (2.16), with initial condition $\xi_i(\theta) = \xi_{i,0}(\theta)$ for $\theta \in [-\delta, 0)$ ($i = 1, ..., N$), $(z(0), w(0)) = (z_0, w_0)$, $x(\theta) = x_0(\theta)$ for $\theta \in [-r, 0]$, $u(\theta) = u_0(\theta)$ for $\theta \in [-r-\tau, 0)$ corresponding to input $e \in L^{\infty}_{loc}\left(\Re_+; \Re^k\right)$.

We prove next that the solution exists for all $t \geq 0$. In order to prove that the solution exists for all $t \geq 0$, it suffices to show that the solution exists and is bounded for all $t \in [0, \tau_1)$. Indeed, if the solution exists and is bounded for all $t \in [0, \tau_1)$ then $x(\tau_1)$ and $z(\tau_1)$ can be uniquely defined and consequently $w(\tau_1)$ can be uniquely defined (by means of (2.14)). Therefore, all arguments can be repeated to the interval $[\tau_1, \tau_2)$ and in the same way we obtain existence of solution for all intervals $[\tau_i, \tau_{i+1})$ ($i = 0,1,2,...$).

Standard results in ordinary differential equations guarantee that the system

$$\dot{z}(t) = f(z(t), u(t-r-\tau)) + \hat{k}(z(t), w(t), u(t-r-\tau)) \qquad (3.2)$$
$$\dot{w}(t) = \nabla h(z(t)) f(z(t), u(t-r-\tau))$$

has a local solution defined on $[0, t_1)$ for some $t_1 \in (0, \min(\tau_1, r+\tau)]$. By virtue of (3.1) and Lemma 2.1, it follows that the solution of (3.2) satisfies the following estimate:

$$V(z(t)) \leq \max(V(z_0), b) \qquad (3.3)$$

for all $t \geq 0$ for which the solution of (3.2) exists. Define the non-decreasing function:

$$\Omega(s) := \max\{|\nabla h(z) f(z,u)| : (z,u) \in \Re^n \times U, V(z) \leq s\}, \text{ for all } s \geq \min(V(z) : z \in \Re^n) \qquad (3.4)$$

which is well-defined by virtue of the facts that $U \subseteq \Re^m$ is compact and $V \in C^2(\Re^n; \Re_+)$ is a radially unbounded function. It follows from definition (3.4) and inequality (3.3), that the solution of (3.2) satisfies the following estimate for all $t \in [0, t_1)$:

$$|w(t)| \leq |w_0| + T_s \Omega(\max(V(z_0), b)) \qquad (3.5)$$

A standard contradiction argument shows that the solution of (3.2) exists and satisfies (3.3), (3.5) for all $t \in [0, \min(\tau_1, r+\tau))$.

Next consider the solution of the system (2.15), (2.16). System (2.15), (2.16) is a system described by Integral Delay Equations (IDEs) with input $z(t)$. The existence of $t_1 \in (0, \min(\tau_1, r+\tau)]$ for which the solution of system (2.15), (2.16) is uniquely defined on $[0, t_1)$ is a direct consequence of Theorem 2.1 in [12] (in conjunction with the fact that assumptions (H1), (H2) in [12] hold for system (2.15), (2.16)). Using the fact that the inequality $q\left(\frac{|\xi|}{\psi(z)}\right)\xi \leq K\psi(z)$ holds for all $(\xi, z) \in \Re^n \times \Re^n$ in conjunction with definition $p(s) := \max\{|f(z,u)| : (z,u) \in \Re^n \times U, |z| \leq s\}$, we obtain the estimate

$$|\xi_j(t)| \leq K\psi(z(t)) + \delta p(K\psi(z(t))), \quad j = 1,...,N \qquad (3.6)$$

for $t \in [0, t_1)$ a.e.. The fact that system (2.15), (2.16) satisfies the Boundedness-Implies-Continuation property (a consequence of Theorem 2.1 in [12]) in conjunction with estimates (3.3), (3.6) shows that the solution of (2.15), (2.16) exists and satisfies (3.3) and (3.6) for $t \in [0, \min(\tau_1, r+\tau))$ a.e..

Finally, the solution of (1.1) exists locally and by virtue of (2.1) and Lemma 2.1 satisfies the estimate:

$$V(x(t)) \leq \max(V(x_0(0)), R) \qquad (3.7)$$

for all $t \geq 0$ for which the solution of (1.1) exists. A standard contradiction argument in conjunction with the fact that $V \in C^2(\Re^n; \Re_+)$ is a radially unbounded function guarantees that the solution of (1.1) exists and satisfies (3.7) for all $t \in [0, \min(\tau_1, r + \tau))$.

If $r + \tau < \tau_1$ then all arguments may be repeated for the interval $t \in [r + \tau, \min(\tau_1, 2r + 2\tau))$ and continuing in this way we show that the solution of (1.1), (2.12), (2.13), (2.14), (2.15), (2.16) exists for all $t \in [0, \tau_1)$.

Lemma 2.1 in conjunction with (2.1) and (3.1) implies there exists $T \in C^0(\Re^n; \Re_+)$ such that (3.3), (3.7) hold for all $t \geq 0$ and

$$V(x(t)) \leq R \text{ for all } t \geq T(x_0(0)) \text{ and } V(z(t)) \leq b \text{ for all } t \geq T(z_0) \qquad (3.8)$$

Indeed, the above conclusions for $V(x(t))$ are direct consequences of Lemma 2.1. The above conclusions for $V(z(t))$ are consequences of Lemma 2.1 applied to system (2.12) with $(w, u)$ as inputs. Inequalities (3.6), (3.8) show that

$$x(t) \in S_1, \ z(t) \in S_2 \ \xi_j(t) \in S_3 \ (j = 1, \ldots, N), \text{ for } t \geq \max(T(x_0(0)), T(z_0)) \text{ a.e.} \qquad (3.9)$$

where $S_1 := \{x \in \Re^n : V(x) \leq R\}$, $S_2 := \{x \in \Re^n : V(x) \leq b\}$, $S_3 := \{x \in \Re^n : |x| \leq K\psi(z) + \delta p(K\psi(z)), z \in S_2\}$.

Equation (1.1) implies that
$x(t - r + j\delta) = x(t - r + (j - 1)\delta) + \int_0^\delta f(x(t - r + s + (j - 1)\delta), u(t + (j - 1)\delta - \tau - r + s)) ds$ for $j = 1, \ldots, N$ and $t \geq r$.
Using the previous equation in conjunction with (2.1), implication (2.10) and Lemma 2.1 (which imply that $q\left(\frac{|x(t - r + s + (j - 1)\delta)|}{\psi(x(t - r))}\right) = 1$ for all $s \geq 0$, $t \geq r$ and $j = 1, \ldots, N$), we get from (2.15):

$$\begin{aligned}
\xi_j(t) - x(t - r + j\delta) &= q\left(\frac{|\xi_{j-1}(t)|}{\psi(z(t))}\right)\xi_{j-1}(t) - q\left(\frac{|x(t - r + (j - 1)\delta)|}{\psi(x(t - r))}\right)x(t - r + (j - 1)\delta) \\
&+ \int_0^\delta f\left(q\left(\frac{|\xi_j(t + s - \delta)|}{\psi(z(t))}\right)\xi_j(t + s - \delta), u(t + (j - 1)\delta - \tau - r + s)\right)ds \\
&- \int_0^\delta f\left(q\left(\frac{|x(t - r + s + (j - 1)\delta)|}{\psi(x(t - r))}\right)x(t - r + s + (j - 1)\delta), u(t + (j - 1)\delta - \tau - r + s)\right)ds
\end{aligned} \qquad (3.10)$$

for all $j = 1, \ldots, N$ and $t \geq r$. Equation (3.10) in conjunction with (3.9) the fact that the inequality $q\left(\frac{|\xi|}{\psi(z)}\right)\xi \leq K\psi(z)$ holds for all $(\xi, z) \in \Re^n \times \Re^n$ and the definitions

$$M_1^q := \sup\left\{\frac{\left|q\left(\frac{|\xi|}{\psi(z)}\right)\xi - q\left(\frac{|x|}{\psi(z)}\right)x\right|}{|x-\xi|} : x \in S_1, \xi \in S_3, z \in S_2, x \neq \xi\right\},\quad M_2^q := \sup\left\{\frac{\left|q\left(\frac{|x|}{\psi(z)}\right)x - q\left(\frac{|x|}{\psi(w)}\right)x\right|}{|z-w|} : x, w \in S_1, z \in S_2, z \neq w\right\},$$

$$M_1^f := \sup\left\{\frac{|f(x,u) - f(\xi,u)|}{|x-\xi|} : x \in S_1, \xi \in S_4, u \in U, x \neq \xi\right\}, \quad S_4 := \{x \in \Re^n : |x| \leq K\psi(z), z \in S_2\}, \text{ implies that:}$$

$$\begin{aligned}|\xi_j(t) - x(t - r + j\delta)| &\leq M_1^q |\xi_{j-1}(t) - x(t - r + (j-1)\delta)| + (1 + \delta M_1^f) M_2^q |z(t) - x(t - r)| \\ &+ \delta M_1^f M_1^q \sup_{-\delta \leq \theta \leq 0}(|\xi_j(t + \theta) - x(t + \theta - r + j\delta)|)\end{aligned} \tag{3.11}$$

for all $j = 1, \ldots, N$ and $t \geq \max(T(x_0(0)), T(z_0)) + \max(r, \delta)$. Using (3.11) in conjunction with (2.11) we obtain:

$$\begin{aligned}\sup_{t_0 \leq t \leq T}(|\xi_j(t) - x(t - r + j\delta)|e^{\sigma t}) &\leq \frac{M_1^q}{1 - \delta M_1^q M_1^f e^{\sigma\delta}} \sup_{t_0 \leq t \leq T}(|\xi_{j-1}(t) - x(t - r + (j-1)\delta)|e^{\sigma t}) \\ &+ \frac{(1 + \delta M_1^f) M_2^q}{1 - \delta M_1^q M_1^f e^{\sigma\delta}} \sup_{t_0 \leq t \leq T}(|z(t) - x(t - r)|e^{\sigma t}) + \sup_{t_0 - \delta \leq t \leq t_0}(|\xi_j(t) - x(t - r + j\delta)|e^{\sigma t})\end{aligned} \tag{3.12}$$

for all $j = 1, \ldots, N$ and $T \geq t_0 := \max(T(x_0(0)), T(z_0)) + \max(r, \delta)$. Using (3.12) we conclude that the following estimate holds

$$\begin{aligned}\sup_{t_0 \leq s \leq t}(|\xi_j(s) - x(s + \tau)|e^{\sigma s}) &\leq \left(\frac{(1 + \delta M_1^f) M_2^q}{1 - \delta M_1^q M_1^f e^{\sigma\delta}}(1 + \ldots + \lambda^{j-1}) + \lambda^j\right) \sup_{t_0 \leq s \leq t}(|z(s) - x(s - r)|e^{\sigma s}) \\ &+ \sum_{l=1}^{j} \lambda^{j-l} \sup_{t_0 - \delta \leq s \leq t_0}(|\xi_l(s) - x(s - r + l\delta)|e^{\sigma s})\end{aligned} \tag{3.13}$$

for all $j = 1, \ldots, N$ and $t \geq t_0 := \max(T(x_0(0)), T(z_0)) + \max(r, \delta)$ with $\lambda := \frac{M_1^q}{1 - \delta M_1^q M_1^f e^{\sigma\delta}}$.

Next consider the evolution of the mapping $t \to P(x(t))$. Inequality (2.2) and (3.9) imply that the following differential inequality holds for $t \geq \tau + \max(T(x_0(0)), T(z_0))$ a.e.:

$$\frac{d}{dt}(P(x(t))) \leq -2\mu|x(t)|^2 + \sqrt{n} M_2^f \widetilde{P}|x(t)||x(t) - \xi_N(t - \tau)| \tag{3.14}$$

where $M_2^f := \sup\left\{\frac{|f(x, k(x)) - f(x, k(\xi))|}{|x - \xi|} : x \in S_1, \xi \in S_3, x \neq \xi\right\}$, $\widetilde{P} := \max\{|\nabla^2 P(x)| : x \in co(S_1)\}$ and $co(S_1)$ denotes the convex hull of $S_1$. Since $P \in \Re^{n \times n}$ is a positive definite matrix. Completing the squares, integrating and noticing that there exists a constant $K_1 > 0$ with $K_1|x|^2 \leq x'Px \leq \frac{\sqrt{n}}{2}\widetilde{P}|x|^2$ for all $x \in S_1$, we obtain the following estimate for $t \geq \tau + t_0$, $t_0 := \max(T(x_0(0)), T(z_0)) + \max(r, \delta)$:

$$|x(t)| \leq e^{-\frac{2\mu}{\sqrt{n\widetilde{P}}}(t-t_0-\tau)} \sqrt{\frac{\sqrt{n\widetilde{P}}}{2K_1}} |x(t_0+\tau)| + \sqrt{\frac{n\sqrt{n\widetilde{P}}}{K_1}} \frac{M_2^f \widetilde{P}}{2\mu} \sup_{t_0+\tau \leq s \leq t} \left( e^{-\frac{\mu}{\sqrt{n\widetilde{P}}}(t-s)} |x(s) - \xi_N(s-\tau)| \right) \quad (3.15)$$

Since $\sigma \leq \frac{\mu}{\sqrt{n\widetilde{P}}}$ (see (2.11)), we obtain from (3.15) and (3.13) for $t \geq \tau + t_0$:

$$\sup_{t_0+\tau \leq s \leq t} \left( |x(s)| e^{\sigma s} \right) \leq \sqrt{\frac{\sqrt{n\widetilde{P}}}{2K_1}} |x(t_0+\tau)| e^{\sigma(\tau+t_0)}$$
$$+ \sqrt{\frac{n\sqrt{n\widetilde{P}}}{K_1}} \frac{M_2^f \widetilde{P}}{2\mu} \left( \frac{(1+\delta M_1^f) M_2^q}{1-\delta M_1^q M_1^f e^{\sigma\delta}} \left(1 + \ldots + \lambda^{N-1}\right) + \lambda^N \right) \sup_{t_0 \leq s \leq t-\tau} \left( |z(s) - x(s-r)| e^{\sigma s} \right) \quad (3.16)$$
$$+ \sqrt{\frac{n\sqrt{n\widetilde{P}}}{K_1}} \frac{M_2^f \widetilde{P}}{2\mu} \sum_{l=1}^{N} \lambda^{N-l} \sup_{t_0-\delta \leq s \leq t_0} \left( |\xi_l(s) - x(s-r+l\delta)| e^{\sigma s} \right)$$

Next we establish the following inequality:

$$(z-x)' Q\left(f(z,u) + \hat{k}(z,h(x),u) - f(x,u)\right) \leq -c\omega |z-x|^2, \text{ for all } (x,z,u) \in S_1 \times S_2 \times U \quad (3.17)$$

Notice that inequality (2.3) and definitions (2.7), (2.8), (2.9) imply that (3.17) holds for the case $V(z) \leq a$. Therefore, we focus on the case $a < V(z) \leq b$. Definition (2.8) gives:

$$(z-x)' Q\left(f(z,u) + \hat{k}(z,h(x),u) - f(x,u)\right)$$
$$\leq (z-x)' Q\left(f(z,u) + L(h(z) - h(x)) - f(x,u)\right) - \frac{\varphi(z,h(x),u)}{|\nabla V(z)|^2} \nabla V(z) Q(z-x) \quad (3.18)$$

Inequalities (2.3), (3.18) and the fact that $\varphi(z,h(x),u) \geq 0$ implies that (3.17) holds if $\nabla V(z) Q(z-x) \geq 0$. Moreover, inequalities (2.3), (3.18) show that (3.17) holds if $\varphi(z,h(x),u) = 0$. It remains to consider the case $\nabla V(z) Q(z-x) < 0$ and $\varphi(z,h(x),u) > 0$. In this case, definition (2.9) implies $\varphi(z,h(x),u) = \nabla V(z) f(z,u) + W(z) + p(V(z)) \nabla V(z) L(h(z) - h(x)) > 0$. Then, inequality (2.6) gives:

$$\varphi(z,h(x),u)) = \nabla V(z) f(z,u) + p(V(z)) \nabla V(z) L(h(z) - h(x)) + W(z) \leq$$
$$+ (1 - p(V(z))) \nabla V(z) f(z,u) + (1 - p(V(z))) W(z) \quad (3.19)$$
$$+ (1-c) |\nabla V(z)|^2 p(V(z)) \frac{(z-x)' Q\left(f(z,u) + L(h(z) - h(x)) - f(x,u)\right)}{\nabla V(\xi) Q(z-x)}$$

Using (3.19), (2.1) and the fact that $0 \leq p(V(z)) \leq 1$, we obtain:

$$-\frac{\varphi(z,h(x),u)) \nabla V(z) Q(z-x)}{|\nabla V(z)|^2} \leq$$
$$-\frac{1-p(V(z))}{|\nabla V(z)|^2} \nabla V(z) Q(z-x) \left(\nabla V(z) f(z,u) + W(z)\right)$$
$$-(1-c) p(V(z)) (z-x)' Q\left(f(z,u) + L(h(z) - h(x)) - f(x,u)\right)$$
$$\leq -(1-c)(z-x)' Q\left(f(z,u) + L(h(z) - h(x)) - f(x,u)\right)$$

Combining (2.3), (3.18) and the above inequality, we conclude that (3.17) holds.

Next consider the evolution of the mapping $t \to (z(t)-x(t-r))'Q(z(t)-x(t-r))$. Inequality (3.17) and (3.9) imply that the following differential inequality holds for $t \geq r + \max(T(x_0(0)), T(z_0))$ a.e.:

$$\frac{d}{dt}\left((z(t)-x(t-r))'Q(z(t)-x(t-r))\right) \leq -2c\omega|z(t)-x(t-r)|^2 + 2G_2|Q||z(t)-x(t-r)||w(t)-h(x(t-r))| \quad (3.20)$$

where $G_2 := \sup\left\{\frac{|\hat{k}(z,y,u)-\hat{k}(z,w,u)|}{|y-w|} : y, w \in \Re^k, z \in S_2, u \in U, y \neq w\right\}$. Since $Q \in \Re^{n \times n}$ is a positive definite matrix there exists a constant $0 < K_2 \leq |Q|$ with $K_2|x|^2 \leq x'Qx$ for all $x \in \Re^n$. Completing the squares and integrating we obtain the following estimate for $t \geq t_0$, $t_0 := \max(T(x_0(0)), T(z_0)) + \max(r, \delta)$:

$$|z(t)-x(t-r)| \leq e^{-\frac{c\omega}{2|Q|}(t-t_0)}\sqrt{\frac{|Q|}{K_2}}|z(t_0)-x(t_0-r)| + \sqrt{\frac{2|Q|}{K_2}}\frac{G_2|Q|}{c\omega}\sup_{t_0 \leq s \leq t}\left(e^{-\frac{c\omega}{4|Q|}(t-s)}|w(s)-h(x(s-r))|\right) \quad (3.21)$$

Since $\sigma \leq \frac{c\omega}{4|Q|}$ (see (2.11)), we obtain from (3.21) for $t \geq t_0$:

$$\sup_{t_0 \leq s \leq t}\left(e^{\sigma s}|z(s)-x(s-r)|\right) \leq \sqrt{\frac{|Q|}{K_2}}e^{\sigma t_0}|z(t_0)-x(t_0-r)| + \sqrt{\frac{2|Q|}{K_2}}\frac{G_2|Q|}{c\omega}\sup_{t_0 \leq s \leq t}\left(e^{\sigma s}|w(s)-h(x(s-r))|\right) \quad (3.22)$$

Finally, notice that since $\sup_{i \geq 0}(\tau_{i+1}-\tau_i) \leq T_s$ the following estimate holds for every $t \in [\tau_i, \tau_{i+1})$ with $\tau_i \geq t_0$:

$$|w(t)-h(x(t-r))| \leq \sup_{0 \leq s \leq t}|e(s)| + T_s G_1 \sup_{\tau_i \leq s \leq t}|z(s)-x(s-r)| \quad (3.23)$$

where $G_1 := \sup\left\{\frac{|\nabla h(x)f(x,u)-\nabla h(z)f(z,u)|}{|x-z|} : x \in S_1, z \in S_2, u \in U, x \neq z\right\}$. Notice that from the inequalities $t \leq \tau_i + T_s$, $\tau_i \leq t_0 + T_s$ and (3.23) we obtain for all $t \geq t_0 + T_s$:

$$\sup_{t_0+T_s \leq s \leq t}\left(e^{\sigma s}|w(s)-h(x(s-r))|\right) \leq e^{\sigma t}\sup_{0 \leq s \leq t}|e(s)| + T_s G_1 e^{\sigma T_s}\sup_{t_0 \leq s \leq t}\left(e^{\sigma s}|z(s)-x(s-r)|\right) \quad (3.24)$$

Combining (3.22) and (3.24) we get for all $t \geq t_0 + T_s$:

$$\sup_{t_0 \leq s \leq t}\left(e^{\sigma s}|z(s)-x(s-r)|\right) \leq \sqrt{\frac{|Q|}{K_2}}e^{\sigma t_0}|z(t_0)-x(t_0-r)| + \sqrt{\frac{2|Q|}{K_2}}\frac{G_2|Q|}{c\omega}e^{\sigma t}\sup_{0 \leq s \leq t}(|e(s)|)$$
$$+ T_s G_1 e^{\sigma T_s}\sqrt{\frac{2|Q|}{K_2}}\frac{G_2|Q|}{c\omega}\sup_{t_0 \leq s \leq t}\left(e^{\sigma s}|z(s)-x(s-r)|\right) + \sqrt{\frac{2|Q|}{K_2}}\frac{G_2|Q|}{c\omega}\sup_{t_0 \leq s \leq t_0+T_s}\left(e^{\sigma s}|w(s)-h(x(s-r))|\right) \quad (3.25)$$

It follows from (2.11), (3.25) that the following estimate holds for all $t \geq t_0 + T_s$:

$$\sup_{t_0 \leq s \leq t}\left(e^{\sigma s}|z(s)-x(s-r)|\right) \leq \frac{c\omega\sqrt{|Q|}}{c\omega\sqrt{K_2}-T_s G_1 G_2 |Q| e^{\sigma T_s}\sqrt{2|Q|}} e^{\sigma t_0}|z(t_0)-x(t_0-r)|$$
$$+\frac{G_2|Q|\sqrt{2|Q|}}{c\omega\sqrt{K_2}-T_s G_1 G_2 |Q| e^{\sigma T_s}\sqrt{2|Q|}}\left(e^{\sigma t}\sup_{0\leq s\leq t}(|e(s)|)+\sup_{t_0\leq s\leq t_0+T_s}\left(e^{\sigma s}|w(s)-h(x(s-r))|\right)\right)$$

(3.26)

Combining (3.16) and (3.26), we obtain the following inequality for all $t \geq t_0 + T_s$:

$$|x(t)|e^{\sigma t} \leq A_1|x(t_0+\tau)|e^{\sigma(\tau+t_0)} + A_2 \sup_{t_0 \leq s \leq t_0+T_s}\left(|w(s)-h(x(s-r))|e^{\sigma s}\right)$$
$$+ A_3 e^{\sigma t_0}|z(t_0)-x(t_0-r)| + A_4 e^{\sigma t}\sup_{0\leq s\leq t}(|e(s)|) + A_5 \sum_{l=1}^{N}\sup_{t_0-\delta\leq s\leq t_0}\left(|\xi_l(s)-x(s-r+l\delta)|e^{\sigma s}\right)$$

(3.27)

for appropriate constants $A_i > 0$ ($i = 1,\ldots,5$). Combining (3.13), (3.24), (3.26), (3.27), using (2.16) and the fact that $k: \Re^n \to U$ is globally Lipschitz with $k(0) = 0$ and defining

$$T_0 := \max(T(x_0(0)), T(z_0)) + r + 2\tau + T_s$$

(3.28)

we obtain the following estimate for all $t \geq 0$:

$$\sup_{-r\leq s\leq 0}(|x(t+s)|)+|w(t)|+|z(t)|+\sum_{j=1}^{N}\sup_{-\delta\leq s<0}(|\xi_j(t+s)|)+\sup_{-r-\tau\leq s<0}(|u(t+s)|) \leq$$
$$Ae^{-\sigma(t-T_0)}\left(\sup_{-r\leq s\leq T_0}(|x(s)|)+\sup_{0\leq s\leq T_0}(|z(s)|)+\sup_{0\leq s\leq T_0}(|w(s)|)+\sum_{j=1}^{N}\sup_{-\delta\leq s<T_0}(|\xi_j(s)|)+\sup_{-r-\tau\leq s<T_0}(|u(s)|)\right)+\gamma\sup_{0\leq s\leq t}(|e(s)|)$$

(3.29)

for appropriate constants $A, \gamma > 0$.

Estimates (3.3), (3.6), (3.7) and the fact that $k: \Re^n \to U$ is globally Lipschitz with $k(0) = 0$ guarantee that there exists a non-decreasing, smooth function $\Delta: \Re_+ \to \Re_+$ such that

$$\sup_{-r\leq s\leq 0}(|x(t+s)|)+|z(t)|+\sum_{j=1}^{N}\sup_{-\delta\leq s<0}(|\xi_j(t+s)|)+\sup_{-r-\tau\leq s<0}(|u(t+s)|) \leq$$
$$\Delta\left(\sup_{-r\leq s\leq 0}(|x(s)|)+|z(0)|+\sum_{j=1}^{N}\sup_{-\delta\leq s<0}(|\xi_j(s)|)+\sup_{-r-\tau\leq s<0}(|u(s)|)\right)$$

(3.30)

for all $t \geq 0$. Definitions (2.7), (2.8), (2.9) in conjunction with the fact that $q(s) \leq K$ for all $s \geq 0$, guarantee the existence of a smooth, non-decreasing function $G: \Re_+ \to \Re_+$ such that the following inequality holds for every $R \geq 0$:

$$\left|f\left(q\left(\frac{|\xi|}{\psi(z)}\right)\xi, u_1\right)\right|+|f(x,u_2)|+|\nabla h(z)f(z,u_3)|+|f(z,u_3)+\hat{k}(z,w,u_3)|+|h(\zeta)|,$$
$$\leq G(R)(|\xi|+|x|+|z|+|w|+|u_1|+|u_2|+|u_3|+|\zeta|)$$

for all $(\xi,\zeta,x,z) \in (\Re^n)^4$, $w \in \Re^k$, $(u_1,u_2,u_3) \in U^3$ with $|\xi|+\max(|x|,|\zeta|)+|z|+\max_{i=1,2,3}|u_i| \leq R$ (3.31)

In order to finish the proof, we notice that is suffices to prove that there exist smooth functions $p_j : \Re_+ \times \Re_+ \to \Re_+$ for which $p_j(R,\cdot)$ and $p_j(\cdot,t)$ are non-decreasing for every fixed $(R,t) \in \Re_+ \times \Re_+$ ($j = 1,2,3,4$) such that the following estimates

$$|x(t)| \leq p_1(R,t) \left( \sup_{-r \leq s \leq 0} (|x(s)|) + |z(0)| + |w(0)| + \sum_{j=1}^{N} \sup_{-\delta \leq s < 0} (|\xi_j(s)|) + \sup_{-r-\tau \leq s < 0} (|u(s)|) + \sup_{0 \leq s \leq t} (|e(s)|) \right) \quad (3.32)$$

$$|z(t)| \leq p_2(R,t) \left( \sup_{-r \leq s \leq 0} (|x(s)|) + |z(0)| + |w(0)| + \sum_{j=1}^{N} \sup_{-\delta \leq s < 0} (|\xi_j(s)|) + \sup_{-r-\tau \leq s < 0} (|u(s)|) + \sup_{0 \leq s \leq t} (|e(s)|) \right) \quad (3.33)$$

$$|w(t)| \leq p_3(R,t) \left( \sup_{-r \leq s \leq 0} (|x(s)|) + |z(0)| + |w(0)| + \sum_{j=1}^{N} \sup_{-\delta \leq s < 0} (|\xi_j(s)|) + \sup_{-r-\tau \leq s < 0} (|u(s)|) + \sup_{0 \leq s \leq t} (|e(s)|) \right) \quad (3.34)$$

$$|u(t)| + \sum_{j=1}^{N} |\xi_j(t)| \leq p_4(R,t) \left( \sup_{-r \leq s \leq 0} (|x(s)|) + |z(0)| + |w(0)| + \sum_{j=1}^{N} \sup_{-\delta \leq s < 0} (|\xi_j(s)|) + \sup_{-r-\tau \leq s < 0} (|u(s)|) + \sup_{0 \leq s \leq t} (|e(s)|) \right) \quad (3.35)$$

hold for all $t \geq 0$ with $R := \sup_{-r \leq s \leq 0} (|x(s)|) + |z(0)| + \sum_{j=1}^{N} \sup_{-\delta \leq s < 0} (|\xi_j(s)|) + \sup_{-r-\tau \leq s < 0} (|u(s)|)$. Indeed, if estimates (3.32), (3.33), (3.34), (3.35) hold then by virtue of (3.28), (3.29) we conclude that (2.17) holds with $\Gamma := \gamma$ and $C(s) := sAe^{\sigma \widetilde{T}(s)} \max_{j=1,2,3,4} (p_j(s, \widetilde{T}(s)))$ for all $s \geq 0$, where $\widetilde{T} : \Re_+ \to \Re_+$ is a smooth, non-decreasing function that satisfies

$$r + 2\tau + T_s + \max\{T(x) : |x| \leq s\} \leq \widetilde{T}(s), \text{ for all } s \geq 0 \quad (3.36)$$

For convenience we use $G$ in order to denote $G(\Delta(R))$ with $R := \sup_{-r \leq s \leq 0} (|x(s)|) + |z(0)| + \sum_{j=1}^{N} \sup_{-\delta \leq s < 0} (|\xi_j(s)|) + \sup_{-r-\tau \leq s < 0} (|u(s)|)$. Using (2.12), (3.30), (3.31) we get

$$|z(t)| \leq |z(0)| + G\int_0^t |z(s)|ds + G\int_0^t |w(s)|ds + G\int_0^t |u(s-r-\tau)|ds \quad (3.37)$$

for all $t \geq 0$. Applying the Gronwall-Belman Lemma to (3.37) we obtain

$$|z(t)| \leq e^{2Gt}\left(|z(0)| + G\int_0^t |w(s)|ds + G\int_0^t |u(s-r-\tau)|ds\right) \quad (3.38)$$

for all $t \geq 0$. Using (2.13), (2.14), (3.30), (3.31) we get

$$|w(t)| \leq |w(0)| + G \sup_{-r \leq s \leq t-r} |x(s)| + G \sup_{0 \leq s \leq t} |e(s)| + G\int_{\tau_i}^t |z(s)|ds + G\int_{\tau_i}^t |u(s-r-\tau)|ds \quad (3.39)$$

for all $t \in [\tau_i, \tau_{i+1})$. Since $t \leq \tau_i + T_s$ we obtain from (3.38), (3.39) for all $t \in [\tau_i, \tau_{i+1})$:

$$|w(t)| \leq |w(0)| + G \sup_{-r \leq s \leq t-r} |x(s)| + G \sup_{0 \leq s \leq t} |e(s)| + GT_s e^{2Gt} |z(0)| + G^2 T_s e^{2Gt} \int_0^t |w(s)| ds + G\left(1 + GT_s e^{2Gt}\right) \int_0^t |u(s-r-\tau)| ds$$
(3.40)

Notice that (3.40) holds for all $t \geq 0$. Applying the Gronwall-Belman Lemma to (3.40) and noticing that the Gronwall-Belman Lemma holds not only for continuous functions but also for piecewise continuous functions (and that the mapping $t \to |w(t)|$ is piecewise continuous on $\Re_+$), we get

$$|w(t)| \leq (1+G)\left(1 + GT_s e^{2Gt}\right) e^{G^2 T_s e^{2Gt} t} \left(|w(0)| + \sup_{-r \leq s \leq t-r} |x(s)| + \sup_{0 \leq s \leq t} |e(s)| + |z(0)| + \int_0^t |u(s-r-\tau)| ds\right)$$
(3.41)

for all $t \geq 0$. Combining (3.38) with (3.41) we obtain

$$|w(t)| + |z(t)| \leq a_1(R,t)\left(|w(0)| + \sup_{-r \leq s \leq t-r} |x(s)| + \sup_{0 \leq s \leq t} |e(s)| + |z(0)| + \int_{-r-\tau}^t |u(s)| ds\right)$$
(3.42)

for all $t \geq 0$ with $a_1(R,t) := e^{3Gt}(1+G)\left(2 + GT_s e^{2Gt}\right) e^{G^2 T_s e^{2Gt} t}$, where $G := G(\Delta(R))$.

We next continue with $|\xi_j(t)|$ ($j=1,\ldots,N$). We notice that (2.15) shows the mappings $t \to \xi_j(t)$ are continuous for all $t \geq 0$ ($j=1,\ldots,N$). Using (2.15), (3.30), (3.11), in conjunction with the fact that $q(s) \leq K$ for all $s \geq 0$, we get

$$|\xi_j(t)| \leq K|\xi_{j-1}(t)| + G \int_{t-\delta}^t |\xi_j(s)| ds + G \int_{-r-\tau}^t |u(s)| ds$$
(3.43)

for all $t \geq 0$ and $j=1,\ldots,N$. Inequality (3.43) implies that the following inequality

$$|\xi_j(t)| \leq K|\xi_{j-1}(t)| + G \int_0^t |\xi_j(s)| ds + G \int_{-\delta}^0 |\xi_j(s)| ds + G \int_{-r-\tau}^t |u(s)| ds$$
(3.44)

holds for all $t \geq 0$ and $j=1,\ldots,N$. Applying the Gronwall-Belman Lemma to (3.44) we obtain

$$|\xi_j(t)| \leq (1 + Gte^{Gt})\left(K \sup_{0 \leq s \leq t} \left(|\xi_{j-1}(s)|\right) + G \int_{-r-\tau}^t |u(s)| ds + G\delta \sup_{-\delta \leq s < 0} \left(|\xi_j(s)|\right)\right)$$
(3.45)

for all $t \geq 0$ and $j=1,\ldots,N$. Using induction, inequalities (3.42), (3.45), the fact that $\xi_0(t) = z(t)$ and the fact that $K \geq 1$ we obtain

$$|\xi_j(t)| \leq$$
$$(1 + Gte^{Gt})^j K^j a_1(R,t)\left(|z(0)| + \sup_{-r \leq s \leq t-r} |x(s)| + \sup_{0 \leq s \leq t} |e(s)| + |w(0)| + \left(1 + \frac{2^{j-1}G}{K}\right) \int_{-r-\tau}^t |u(s)| ds + \frac{G\delta}{K} \sum_{l=1}^j \sup_{-\delta \leq s < 0} \left(|\xi_l(s)|\right)\right)$$
(3.46)

for all $t \geq 0$ and $j=1,\ldots,N$. Let $\Lambda \geq 0$ be a constant for which the inequality $|k(x)| \leq \Lambda |x|$ for all $x \in \Re^n$ (since $k: \Re^n \to U$ is globally Lipschitz with $k(0) = 0$ the existence of $\Lambda \geq 0$ is guaranteed). Using (3.46) with $j = N$ and (2.16) we get

$$|u(t)| \leq a_2(R,t)\left(|z(0)| + \sup_{-r \leq s \leq t-r}|x(s)| + \sup_{0 \leq s \leq t}|e(s)| + |w(0)| + \left(1 + \frac{2^{N-1}G}{K}\right)\int_{-r-\tau}^{t}|u(s)|ds + \frac{G\delta}{K}\sum_{l=1}^{N}\sup_{-\delta \leq s < 0}(|\xi_l(s)|)\right) \quad (3.47)$$

for all $t \geq 0$ with $a_2(R,t) := \Lambda(1+Gte^{Gt})^N K^N a_1(R,t)$, where $G := G(\Delta(R))$. Applying the Gronwall-Belman Lemma to (3.44) we obtain

$$|u(t)| \leq a_3(R,t)\left(|z(0)| + \sup_{-r \leq s \leq t-r}|x(s)| + \sup_{0 \leq s \leq t}|e(s)| + |w(0)| + \left(1 + \frac{2^{N-1}G}{K}\right)\int_{-r-\tau}^{0}|u(s)|ds + \frac{G\delta}{K}\sum_{l=1}^{N}\sup_{-\delta \leq s < 0}(|\xi_l(s)|)\right) \quad (3.48)$$

for all $t \geq 0$ with $a_3(R,t) := a_2(R,t)e^{a_2(R,t)\left(1+\frac{2^{N-1}G}{K}\right)t}$, where $G := G(\Delta(R))$. Combining (3.42), (3.46) with (3.48) we obtain:

$$|w(t)| + |z(t)| \leq a_4(R,t)\left(|w(0)| + \sup_{-r \leq s \leq t-r}|x(s)| + \sup_{0 \leq s \leq t}|e(s)| + |z(0)| + \int_{-r-\tau}^{0}|u(s)|ds + \sum_{l=1}^{N}\sup_{-\delta \leq s < 0}(|\xi_l(s)|)\right) \quad (3.49)$$

$$|u(t)| + \sum_{j=1}^{N}|\xi_j(t)| \leq a_5(R,t)\left(|z(0)| + \sup_{-r \leq s \leq t-r}|x(s)| + \sup_{0 \leq s \leq t}|e(s)| + |w(0)| + \int_{-r-\tau}^{0}|u(s)|ds + \sum_{l=1}^{N}\sup_{-\delta \leq s < 0}(|\xi_l(s)|)\right) \quad (3.50)$$

for all $t \geq 0$ with $a_4(R,t) := \left(2 + \frac{(\delta + 2^{N-1})G}{K}\right)(1+ta_3(R,t))a_1(R,t)$,

$a_5(R,t) := (N+\Lambda)(1+Gte^{Gt})^N K^N a_1(R,t)\left(1 + \frac{(\delta + 2^{N-1})G}{K}\right)^2(1+ta_3(R,t))$, where $G := G(\Delta(R))$.

Finally, using (1.1), (3.30), (3.31) we get

$$|x(t)| \leq |x(v)| + G\int_{v}^{t}|x(s)|ds + G\int_{v}^{t}|u(s-\tau)|ds \quad (3.51)$$

for all $T \geq t \geq v \geq 0$. Combining (3.51) with (3.50) we get

$$|x(t)| \leq G(t-v)(1+a_5(R,T))\sup_{v \leq s \leq t}|x(s)| + (1+G+G(t-v)a_5(R,T))\left(M + \sup_{-r \leq s \leq v}|x(s)|\right) \quad (3.52)$$

for all $T \geq t \geq v \geq 0$, where $M := |x(0)| + |z(0)| + |w(0)| + \sup_{0 \leq s \leq T}|e(s)| + \int_{-r-\tau}^{0}|u(s)|ds + \sum_{j=1}^{N}\sup_{-\delta \leq s < 0}(|\xi_j(s)|)$. Let $\theta > 0$ be a constant with $G\theta(1+a_5(R,T)) \leq 1/2$. It follows from (3.52) that the inequality

$$|x(t)| \leq 2(2+G)\left(M + \sup_{-r \leq s \leq v}|x(s)|\right) \quad (3.53)$$

holds for all $T \geq v \geq 0$ and $t \in [v, \min(v+\theta, T)]$. Applying (3.53) repeatedly we get:

$$\sup_{-r \leq s \leq i\theta}|x(s)| \leq M\left(2(2+G) + ... + 2^i(2+G)^i\right) + 2^i(2+G)^i \sup_{-r \leq s \leq 0}|x(s)| \quad (3.54)$$

for all $T \geq 0$ and all integers $i \geq 0$ with $i\theta \leq T$. Selecting $i \geq 0$ so that $i := 1 + [2GT(1+a_5(R,T))]$ and $\theta = T/i$, we get (3.32) with $p_1(R,t) := (1+2Gt(1+a_5(R,t)))(2(2+G))^{(1+2Gt(1+a_5(R,t)))}$. Exploiting (3.32), (3.49), (3.50) we get inequalities (3.33), (3.34), (3.35) for appropriate smooth functions $p_j : \Re_+ \times \Re_+ \to \Re_+$ for which $p_j(R, \cdot)$ and $p_j(\cdot, t)$ are non-decreasing for every fixed $(R, t) \in \Re_+ \times \Re_+$ ($j = 2,3,4$). The proof is complete. ◁

We next continue with the proof of Theorem 2.4.

**Proof of Theorem 2.4:** We first notice that if the solution of the closed-loop system (2.18), (2.19) with (2.12), (2.13), (2.14), (2.15), (2.16), (2.22), (2.23), (2.24) exists for all $t \geq 0$ then the following equalities hold:

$$\theta(t) = a_1(x(t+\tau)), \text{ for all } t \geq \tau_l \quad (3.55)$$

$$\dot{x}(t) = f(x(t), u(t-\tau)), \text{ for all } t \geq \tau_l + \tau \text{ a.e.} \quad (3.56)$$

where $\tau_l = \min\{\tau_i : i \geq 1, \tau_i \geq r\}$. The above equalities are direct consequences of (2.21) and the fact that the equation $\nabla a_1(x) \tilde{f}(x,v) = g(a_1(x), v)$ holds for all $(x, v) \in \Re^n \times \Re^m$. Therefore, by virtue of Theorem 2.2, there exists a locally Lipschitz function $\tilde{C} \in K_\infty$ such that the solution of the closed-loop system (2.18), (2.19) with (2.12), (2.13), (2.14), (2.15), (2.16), (2.22), (2.23), (2.24) satisfies the following estimate for all $t \geq \tau_l + \tau$:

$$\sup_{-r \leq s \leq 0} (|x(t+s)|) + |w(t)| + |z(t)| + \sum_{j=1}^{N} \sup_{-\delta \leq s < 0} (|\xi_j(t+s)|) + \sup_{-r-\tau \leq s < 0} (|u(t+s)|) \leq$$

$$e^{-\sigma(t-\tau_l-\tau)} \tilde{C} \left( \sup_{-r \leq s \leq 0} (|x(\tau_l+\tau+s)|) + |z(\tau_l+\tau)| + |w(\tau_l+\tau)| + \sum_{j=1}^{N} \sup_{-\delta \leq s < 0} (|\xi_j(\tau_l+\tau+s)|) + \sup_{-r-\tau \leq s < 0} (|u(\tau_l+\tau+s)|) \right) \quad (3.57)$$

Since $a_1 : \Re^n \to \Re^l$ is a locally Lipschitz function with $a_1(0) = 0$, there exists a non-decreasing function $c : \Re_+ \to \Re_+$ such that $|a_1(x)| \leq |x| c(|x|)$, for all $x \in \Re^n$. The previous inequality in conjunction with (3.55) and (3.57) implies the existence of a locally Lipschitz function $\overline{C} \in K_\infty$ such that the solution of the closed-loop system (2.18), (2.19) with (2.12), (2.13), (2.14), (2.15), (2.16), (2.22), (2.23), (2.24) satisfies the following estimate for all $t \geq \tau_l + \tau$:

$$\sup_{-r \leq s \leq 0} (|x(t+s)|) + |w(t)| + |\theta(t)| + |z(t)| + \sum_{j=1}^{N} \sup_{-\delta \leq s < 0} (|\xi_j(t+s)|) + \sup_{-r-\tau \leq s < 0} (|u(t+s)|) \leq$$

$$e^{-\sigma(t-\tau_l-\tau)} \overline{C} \left( \sup_{-r \leq s \leq 0} (|x(\tau_l+\tau+s)|) + |z(\tau_l+\tau)| + |w(\tau_l+\tau)| + \sum_{j=1}^{N} \sup_{-\delta \leq s < 0} (|\xi_j(\tau_l+\tau+s)|) + \sup_{-r-\tau \leq s < 0} (|u(\tau_l+\tau+s)|) \right) \quad (3.58)$$

Inequality (2.25) is a direct consequence of estimate (3.58). Therefore, in order to prove Theorem 2.4 it suffices to show that the solution of the closed-loop system (2.18), (2.19) with (2.12), (2.13), (2.14), (2.15), (2.16), (2.22), (2.23), (2.24) exists for all $t \geq 0$.

The arguments for the proof of the existence of the solution the solution of the closed-loop system (2.18), (2.19) with (2.12), (2.13), (2.14), (2.15), (2.16), (2.22), (2.23), (2.24) are exactly the same as those in the proof of the Theorem 2.2 except that:
- we do not use (3.7) but instead we use the fact that (2.18) is forward complete in conjunction with the results in [1], and
- we use the fact that system (2.20) is forward complete.

The proof is complete. ◁

# 4. Illustrative Examples

This section is devoted to the presentation of two nonlinear control systems which can be stabilized by the results of Theorem 2.2 and Theorem 2.4.

**Example 4.1:** This is an example of a two-dimensional nonlinear control system for which all assumptions (H1), (H2), (H3) and (H4) hold. It follows from Theorem 2.2 that the system can be stabilized globally asymptotically and locally exponentially by means of the ISP-O-P-DFC control scheme. The system is described by the equations

$$\begin{aligned} \dot{x}_1(t) &= gx_1(t) - x_1^3(t) + x_2(t) \\ \dot{x}_2(t) &= -x_2^3(t) + u(t-\tau) \\ x(t) &= (x_1(t), x_2(t)) \in \Re^2, u(t) \in [-4\sqrt{2}, 4\sqrt{2}] \subset \Re \end{aligned} \quad (4.1)$$

where $\tau \geq 0$ is the input delay and $g > 0$ is a constant. The measured output is given by the equation:

$$y(\tau_i) = x_1(\tau_i - r) + e(\tau_i) \quad (4.2)$$

where $\{\tau_i\}_{i=0}^{\infty}$ is the sampling partition (a partition of $\Re_+$), $r \geq 0$ is the measurement delay and the input $e: \Re_+ \to \Re$ is the measurement error. We show next that assumptions (H1), (H2), (H3) and (H4) hold for system (4.1) provided that

$$g \leq \frac{1}{167} \quad (4.3)$$

We notice that system (4.1), (4.2) is a system of the form (1.1), (1.2) with $U = [-4\sqrt{2}, 4\sqrt{2}] \subset \Re$ and $f(x,u) := \begin{bmatrix} gx_1 - x_1^3 + x_2 \\ -x_2^3 + u \end{bmatrix}$. We start by showing that assumption (H1) holds with $V(x) := \frac{1}{2}x_1^2 + \frac{1}{2}x_2^2$, $W(x) := \frac{1}{2}V(x)$ and $R := 4$. Indeed, using the inequalities $x_1 x_2 \leq \frac{1}{2}x_1^2 + \frac{1}{2}x_2^2$, $x_2 u \leq \frac{1}{2}x_2^2 + \frac{1}{2}u^2$, $u^2 \leq 32$ and $2V^2(x) \leq x_1^4 + x_2^4$ (which hold for all $(x,u) \in \Re^2 \times [-4\sqrt{2}, 4\sqrt{2}]$), we obtain:

$$\begin{aligned} \nabla V(x) f(x,u) &= gx_1^2 - x_1^4 + x_1 x_2 - x_2^4 + x_2 u \\ &\leq \left(g + \frac{1}{2}\right)x_1^2 - x_1^4 + x_2^2 - x_2^4 + \frac{1}{2}u^2 \\ &\leq \left(g + \frac{1}{2}\right)x_1^2 + x_2^2 + 16 - 2V^2(x) \end{aligned}$$

for all $(x,u) \in \Re^2 \times [-4\sqrt{2}, 4\sqrt{2}]$. Using (4.3) we get $\left(g + \frac{1}{2}\right)x_1^2 + x_2^2 \leq 2V(x)$, which combined with the above inequality gives:

$$\nabla V(x) f(x,u) \leq 2V(x) + 16 - 2V^2(x) \quad (4.4)$$

Inequality (4.4) shows that inequality (2.1) holds with $W(x) := \frac{1}{2}V(x)$ for all $(x,u) \in \Re^2 \times [-4\sqrt{2}, 4\sqrt{2}]$ with $V(x) \geq R = 4$.

We next show that assumption (H2) holds with $P(x) := \frac{1}{2}x_1^2 + \frac{1}{2}(x_2 + 2gx_1)^2$, $K_1 := \frac{1}{4}$, $\mu := \frac{g}{4}$ and

$$k(x) := -\left(5 - \min\left(5, \max\left(4, \frac{1}{2}x_1^2 + \frac{1}{2}x_2^2\right)\right)\right) \Pr_U\left((4g^2 + 1)x_1 + 3gx_2 + 2g(4g^2 - 1)x_1^3 + 12g^2 x_2 x_1^2 + 6gx_2^2 x_1\right) \quad (4.5)$$

where $\Pr_U$ denotes the projection on the set $U = [-4\sqrt{2}, 4\sqrt{2}] \subset \Re$. Notice that since $V(x) := \frac{1}{2}x_1^2 + \frac{1}{2}x_2^2$ we obtain the following estimate

$$\left|(4g^2 + 1)x_1 + 3gx_2 + 2g(4g^2 - 1)x_1^3 + 12g^2 x_2 x_1^2 + 6gx_2^2 x_1\right|$$
$$\leq (4g^2 + 3g + 1)\sqrt{2V(x)} + 2g\left(|4g^2 - 1| + 6g + 3\right) 2V(x)\sqrt{2V(x)}$$

for all $x \in \Re^2$. Taking into account (4.3) and the above inequality we get:

$$\left|(4g^2 + 1)x_1 + 3gx_2 + 2g(4g^2 - 1)x_1^3 + 12g^2 x_2 x_1^2 + 6gx_2^2 x_1\right|$$
$$\leq 2(100g^2 + 67g + 1)\sqrt{2} \leq 4\sqrt{2}$$

for all $x \in \Re^2$ with $V(x) \leq R = 4$. Consequently, we obtain from (4.5) for all $x \in \Re^2$ with $V(x) \leq R = 4$:

$$k(x) = -\left((4g^2 + 1)x_1 + 3gx_2 + 2g(4g^2 - 1)x_1^3 + 12g^2 x_2 x_1^2 + 6gx_2^2 x_1\right) \quad (4.6)$$

Using (4.6) and the definition $P(x) := \frac{1}{2}x_1^2 + \frac{1}{2}(x_2 + 2gx_1)^2$ we obtain for all $x \in \Re^2$ with $V(x) \leq R = 4$:

$$\nabla P(x)f(x, k(x)) = -gx_1^2 - g(x_2 + 2gx_1)^2 - x_1^4 + (x_2 + 2gx_1)\left(-x_2^3 + k(x) + (4g^2 + 1)x_1 - 2gx_1^3 + 3gx_2\right)$$
$$= -gx_1^2 - g(x_2 + 2gx_1)^2 - x_1^4 - (x_2 + 2gx_1)^4$$
$$+ (x_2 + 2gx_1)\left(k(x) + (4g^2 + 1)x_1 + 3gx_2 + 2g(4g^2 - 1)x_1^3 + 12g^2 x_2 x_1^2 + 6gx_2^2 x_1\right) \quad (4.7)$$
$$= -gx_1^2 - g(x_2 + 2gx_1)^2 - x_1^4 - (x_2 + 2gx_1)^4$$

Finally, the inequality

$$x_1^2 + (x_2 + 2gx_1)^2 \geq (1 + 4g^2)x_1^2 + x_2^2 - 4g|x_1 x_2| \geq (1 - 4g^2)x_1^2 + \frac{1}{2}x_2^2$$

in conjunction with (4.3) implies the inequality

$$x_1^2 + (x_2 + 2gx_1)^2 \geq \frac{1}{2}x_1^2 + \frac{1}{2}x_2^2 \quad (4.8)$$

for all $x \in \Re^2$. Inequality (4.8) in conjunction with (4.7) implies inequality (2.2) with $\mu := \frac{g}{4}$. Moreover, inequality (4.8) in conjunction with the definition $P(x) := \frac{1}{2}x_1^2 + \frac{1}{2}(x_2 + 2gx_1)^2$ implies the inequality $K_1|x|^2 \leq P(x)$ with $K_1 := 1/4$.

We next show that assumption (H3) holds with

$$Q := \frac{1}{2}\begin{bmatrix} 1 & -p \\ -p & 1 \end{bmatrix}, \quad L := -\frac{1}{2(1-p^2)}\begin{bmatrix} 2g + 2p(1-pg) + 4p(11+2\sqrt{7})^2 + p \\ 2gp + 4p^2(11+2\sqrt{7})^2 + p^2 + 2(1-pg) \end{bmatrix}, \quad \omega := \frac{p}{4}, \quad b := 7 \qquad (4.9)$$

where $p \in (0,1)$ is a constant. Indeed, by setting $L_1 = -\frac{2g + 2p(1-pg) + 4p(11+2\sqrt{7})^2 + p}{2(1-p^2)}$,

$L_2 = -\frac{2gp + 4p^2(11+2\sqrt{7})^2 + p^2 + 2(1-pg)}{2(1-p^2)}$, $\tilde{e} = \begin{bmatrix} \tilde{e}_1 \\ \tilde{e}_2 \end{bmatrix} = \begin{bmatrix} z_1 - x_1 \\ z_2 - x_2 \end{bmatrix}$, we get:

$$2(\tilde{e})'Q(f(z,u) - f(x,u) + L\tilde{e}_1) =$$
$$g\tilde{e}_1^2 - \tilde{e}_1(z_1^3 - x_1^3) + \tilde{e}_1\tilde{e}_2 + L_1\tilde{e}_1^2$$
$$- pg\tilde{e}_1\tilde{e}_2 - p\tilde{e}_2^2 + p\tilde{e}_1(z_1^3 - x_1^3) - pL_1\tilde{e}_1\tilde{e}_2$$
$$+ p\tilde{e}_1(z_2^3 - x_2^3) - pL_2\tilde{e}_1^2 - \tilde{e}_2(z_2^3 - x_2^3) + L_2\tilde{e}_1\tilde{e}_2 =$$
$$= (L_1 + g - pL_2)\tilde{e}_1^2 + (1 - pg - pL_1 + L_2)\tilde{e}_1\tilde{e}_2 - p\tilde{e}_2^2$$
$$- (1-p)\tilde{e}_1(z_1^3 - x_1^3) - \tilde{e}_2(z_2^3 - x_2^3) + p\tilde{e}_1\tilde{e}_2(z_2^2 + x_2z_2 + x_2^2)$$

Using the fact that $\tilde{e}_1(z_1^3 - x_1^3) \geq 0$, $\tilde{e}_2(z_2^3 - x_2^3) \geq 0$ and the above inequality, we obtain for all $u \in U$, $z, x \in \Re^n$ with $V(z) \leq 7$ and $V(x) \leq 4$:

$$2(\tilde{e})'Q(f(z,u) - f(x,u) + L\tilde{e}_1) \leq (L_1 + g - pL_2)\tilde{e}_1^2 + (1 - pg - pL_1 + L_2)\tilde{e}_1\tilde{e}_2 - p\tilde{e}_2^2 + 2p|\tilde{e}_1\tilde{e}_2|(11+2\sqrt{7})$$

Completing the squares in the above inequality (i.e., using the inequality $2p|\tilde{e}_1\tilde{e}_2|(11+2\sqrt{7}) \leq \frac{p}{2}\tilde{e}_2^2 + 2p(11+2\sqrt{7})^2\tilde{e}_1^2$), we obtain inequality (2.3) with $\omega := \frac{p}{4}$, $b := 7$.

Finally, we show that assumption (H4) holds with $a := 6$ and arbitrary $c \in (0,1)$, provided that

$$p \leq 1/4 \quad \text{and} \quad p(597 + 176\sqrt{7}) \leq \frac{123}{4\sqrt{7}(\sqrt{7}+2)} - 2 - 2g \qquad (4.10)$$

More specifically, we show next that the following inequality

$$\nabla V(z)(f(z,u) + L(h(z) - h(x))) \leq -W(z) \qquad (4.11)$$

holds $u \in U$, $z, x \in \Re^n$ with $6 < V(z) \leq 7$ and $V(x) \leq 4$. Using the definitions $V(x) := \frac{1}{2}x_1^2 + \frac{1}{2}x_2^2$, $W(x) := \frac{1}{2}V(x)$ and the inequalities $z_1z_2 \leq \frac{1}{2}z_1^2 + \frac{1}{2}z_2^2$, $z_2u \leq \frac{1}{2}z_2^2 + \frac{1}{2}u^2$, $u^2 \leq 32$ and $2V^2(z) \leq z_1^4 + z_2^4$ (which hold for all $(z,u) \in \Re^2 \times [-4\sqrt{2}, 4\sqrt{2}]$), we conclude that inequality (4.11) holds provided that the following (more demanding) inequality

$$\left(g + \frac{1}{2}\right)z_1^2 + z_2^2 + L_1z_1(z_1 - x_1) + L_2z_2(z_1 - x_1) + 16 \leq -\frac{1}{2}V(z) + 2V^2(z) \qquad (4.12)$$

holds $u \in U$, $z, x \in \Re^n$ with $6 < V(z) \leq 7$ and $V(x) \leq 4$. Using (4.3) we get $\left(g + \frac{1}{2}\right)z_1^2 + z_2^2 \leq 2V(z)$, which implies that inequality (4.12) holds provided that the following (more demanding) inequality

$$\frac{5}{2}V(z) + (|L_1 z_1| + |L_2 z_2|)|z_1 - x_1| + 16 \leq 2V^2(z) \qquad (4.13)$$

holds $u \in U$, $z, x \in \Re^n$ with $6 < V(z) \leq 7$ and $V(x) \leq 4$. Notice that for all $z, x \in \Re^n$ with $6 < V(z) \leq 7$ and $V(x) \leq 4$, it holds that $|z_1| \leq \sqrt{14}$, $|z_2| \leq \sqrt{14}$, $|z_1 - x_1| \leq \sqrt{14} + 2\sqrt{2}$. The previous inequalities imply that inequality (4.13) holds provided that the following inequality holds:

$$|L_1| + |L_2| \leq \frac{41}{2\sqrt{7}(\sqrt{7} + 2)} \qquad (4.14)$$

Since $|L_1| + |L_2| = (p+1)\frac{2g + 4p(11 + 2\sqrt{7})^2 + p + 2(1 - pg)}{2(1 - p^2)}$, it follows from (4.3) and (4.10) that inequality (4.14) holds.

Define:

$$\hat{k}(z, y, u) := L(z_1 - y), \text{ for all } (z, y, u) \in \Re^2 \times \Re \times [-4\sqrt{2}, 4\sqrt{2}] \text{ with } z_1^2 + z_2^2 \leq 8 \qquad (4.15)$$

$$\hat{k}(z, y, u) := L(z_1 - y) - \frac{\varphi(z, y, u)}{|z|^2}\begin{bmatrix} z_1 \\ z_2 \end{bmatrix}, \text{ for all } (z, y, u) \in \Re^2 \times \Re \times [-4\sqrt{2}, 4\sqrt{2}] \text{ with } z_1^2 + z_2^2 > 8 \qquad (4.16)$$

where $\varphi : \Re^2 \times \Re \times \Re \to \Re_+$ is defined by

$$\varphi(z, y, u) := \max\left(0, \left(g + \frac{1}{4}\right)z_1^2 - z_1^4 + z_1 z_2 - z_2^4 + z_2 u + \frac{1}{4}z_2^2 + p\left(\frac{z_1^2 + z_2^2}{2}\right)(L_1 z_1 + L_2 z_2)(z_1 - y)\right) \qquad (4.17)$$

and $p : \Re_+ \to [0,1]$ is a locally Lipschitz function that satisfies $p(s) = 1$ for all $s \geq b$ and $p(s) = 0$ for all $s \leq a$. Let $q : \Re \to \Re_+$ be the continuously differentiable function with $q(s) := 2s^{-1} - s^{-2}$ for $s > 1$ and $q(s) = 1$ for $s \leq 1$ with $sq(s) \leq K = 2$ for $s \geq 1$. Let $\psi : \Re^n \to [1, +\infty)$ be a smooth function defined by $\psi(z) := \frac{1 + 2\sqrt{14}}{2} + \frac{1}{2}|z|^2$, which satisfies implication (2.10) with $V(x) := \frac{1}{2}x_1^2 + \frac{1}{2}x_2^2$ and $b = 7$. Computing all constants involved in (2.11) and using all previous definitions, Theorem 2.2 implies that there exist a constant $\Gamma > 0$ and a locally Lipschitz function $C \in K_\infty$ such that for every partition $\{\tau_i\}_{i=0}^\infty$ of $\Re_+$ with $\sup_{i \geq 0}(\tau_{i+1} - \tau_i) \leq T_s$, $v \in L^\infty(\Re_+; \Re)$, $\xi_{i,0} \in L^\infty([-\delta, 0); \Re^2)$ ($i = 1, ..., N$), $(z_0, w_0) \in \Re^2 \times \Re$, $x_0 \in C^0([-r, 0]; \Re^2)$, $u_0 \in L^\infty([-r - \tau, 0); [-4\sqrt{2}, 4\sqrt{2}])$, the solution of (4.1), (4.2) with

$$\dot{z}(t) = \begin{bmatrix} gz_1(t) - z_1^3(t) + z_2(t) \\ -z_2^3(t) + u(t - r - \tau) \end{bmatrix} + \hat{k}(z(t), w(t), u(t - r - \tau)), \text{ for } t \geq 0 \qquad (4.18)$$

$$\dot{w}(t) = gz_1(t) - z_1^3(t) + z_2(t), \text{ for } t \in [\tau_i, \tau_{i+1}), i \geq 0 \qquad (4.19)$$

$$w(\tau_i) = y(\tau_i), \text{ for } i \geq 1 \qquad (4.20)$$

$$\xi_j(t) = q\left(\frac{|\xi_{j-1}(t)|}{\psi(z(t))}\right)\xi_{j-1}(t) + \int_0^\delta f\left(q\left(\frac{|\xi_j(t+s-\delta)|}{\psi(z(t))}\right)\xi_j(t+s-\delta), u(t+(j-1)\delta - \tau - r + s)\right)ds, \text{ for } t \geq 0, j = 1,...,N$$

(4.21)

with $\xi_0(t) = z(t)$ and

$$u(t) = k(\xi_N(t)), \text{ for } t \geq 0 \tag{4.22}$$

initial condition $\xi_j(\theta) = \xi_{j,0}(\theta)$ for $\theta \in [-\delta, 0)$ ($j = 1,...,N$), $(z(0), w(0)) = (z_0, w_0)$, $x(\theta) = x_0(\theta)$ for $\theta \in [-r, 0]$, $u(\theta) = u_0(\theta)$ for $\theta \in [-r-\tau, 0)$, exists and satisfies estimate (2.17) for all $t \geq 0$, provided that $\delta = \frac{r+\tau}{N}$, $N > 0$ is an integer and $T_s > 0$, $\sigma > 0$ are constants so that:

$$\sigma < \min\left(\frac{g}{8(1+2g^2)\sqrt{2}}, \frac{p}{8(1+p)}\right), \quad \delta e^{\sigma\delta} < \frac{27(1+2\sqrt{14})}{(27+54\sqrt{14}+32\sqrt{2})(3034+g)} \text{ and } T_s e^{\sigma T_s} < \frac{\sqrt{2}p\sqrt{1-p}}{12(43+g)|L|(1+p)^{3/2}}$$

(4.23)

Of course, it should be noted that inequalities (4.23) are highly conservative. The control practitioner can use inequalities (4.23) only as a first step for the selection of the parameters. The next step is the determination of the optimal values of the parameters by means of extensive numerical experiments. ◁

**Example 4.2:** Consider the nonlinear system

$$\begin{aligned}\dot{x}_1(t) &= gx_1(t) - x_1^3(t) + x_2(t) \\ \dot{x}_2(t) &= -x_2^3(t) + p(x_3(t)) + v_1(t-\tau) \\ \dot{x}_3(t) &= v_2(t-\tau) \\ x(t) &= (x_1(t), x_2(t), x_3(t)) \in \Re^3, v(t) \in \Re^2\end{aligned} \tag{4.24}$$

where $\tau > 0$ is the input delay, $p \in C^\infty(\Re; \Re)$ is a smooth function with $p(0) = 0$ and $g > 0$ is a constant. The measured output is given by the equation:

$$y(\tau_i) = \begin{bmatrix} y_1(\tau_i) \\ y_2(\tau_i) \end{bmatrix} = \begin{bmatrix} x_1(\tau_i - r) \\ x_3(\tau_i - r) \end{bmatrix} \tag{4.25}$$

where $\{\tau_i\}_{i=0}^\infty$ is the sampling partition (a partition of $\Re_+$) and $r \geq 0$ is the measurement delay. All assumptions of Theorem 2.4 are satisfied with

$$a_1(x) := x_3, \quad l = 1, \quad a_2(\theta) = \begin{bmatrix} -p(\theta) \\ -\theta \end{bmatrix}, \quad g(\theta, v) = v_2, \quad \Phi(y, v) := y_2 + \int_{-r-\tau}^0 v_2(s)ds, \quad U = [-4\sqrt{2}, 4\sqrt{2}] \times \{0\} \tag{4.26}$$

under the assumption that inequality (4.3) holds. Indeed, notice that after an initial transient period system (4.24) with

$$v(t) = u(t) - \begin{bmatrix} p(\theta(t)) \\ \theta(t) \end{bmatrix}, \text{ for } t \geq 0 \tag{4.27}$$

$$\dot{\theta}(t) = -\theta(t) + u_2(t), \text{ for } t \in [\tau_i, \tau_{i+1}), \quad i \geq 0 \tag{4.28}$$

$$\theta(\tau_i) = y_2(\tau_i) + \int_{\tau_i - r - \tau}^{\tau_i} v_2(s)ds, \text{ for } i \geq 1 \tag{4.29}$$

is transformed to system (4.1) with the additional equation:

$$\dot{x}_3(t) = -x_3(t) + u_2(t-\tau) \tag{4.30}$$

It follows from Theorem 2.4 that for every partition $\{\tau_i\}_{i=0}^{\infty}$ of $\Re_+$ with $\sup_{i\geq 0}(\tau_{i+1}-\tau_i) \leq T_s$, the solution of the closed-loop system (4.24) with (4.27), (4.28), (4.29), (4.18), (4.19), (4.20), (4.21) and (4.22) with $\xi_0(t) = z(t)$ satisfies (2.25) provided that $\delta = \frac{r+\tau}{N}$, $N>0$ is an integer and $T_s > 0$, $\sigma > 0$ are constants so that (4.23) holds. A closer inspection of the closed-loop system (4.24) with (4.27), (4.28), (4.29), (4.18), (4.19), (4.20), (4.21), (4.22) and $\xi_0(t) = z(t)$ reveals that not only the exponential attractivity property holds for the closed-loop system but also the properties of Lagrange and Lyapunov stability (see [8]).  ◁

## 5. Concluding Remarks

In this work we have shown that the ISP-O-P-DFC control scheme can be applied to nonlinear systems with a compact absorbing set. The results guarantee the following properties are required in closed loop,
- global asymptotic stability and exponential convergence for the disturbance-free case,
- robustness with respect to perturbations of the sampling schedule,
- robustness with respect to measurement errors,

even when the full state is not measured and when the measurement is sampled and possibly delayed.

More remains to be done for the class of systems which can be transformed to a nonlinear system with a compact absorbing set by means of a preliminary predictor feedback. Although Theorem 2.4 guarantees global exponential attractivity in the absence of measurement errors, additional assumptions must be employed for the global asymptotic stability and exponential convergence in the disturbance-free case.

An extension of the previous results to the case where the control is applied through a zero order hold device is also an open problem and it is under investigation by the authors.

## References

[1] Angeli, D. and E.D. Sontag, "Forward Completeness, Unbounded Observability and their Lyapunov Characterizations", *Systems and Control Letters*, 38(4-5), 1999, 209-217.
[2] Ahmed-Ali, T., I. Karafyllis and F. Lamnabhi-Lagarrigue, "Global Exponential Sampled-Data Observers for Nonlinear Systems with Delayed Measurements", *Systems and Control Letters*, 62(7), 2013, 539-549.
[3] Bekiaris-Liberis, N. and M. Krstic, "Compensation of time-varying input and state delays for nonlinear systems", *Journal of Dynamic Systems, Measurement, and Control*, 134, paper 011009, 2012.
[4] Bekiaris-Liberis, N. and M. Krstic, "Compensation of state-dependent input delay for nonlinear systems", *IEEE Transactions on Automatic Control*, 58, 2013, 275-289.
[5] Bekiaris-Liberis, N. and M. Krstic, "Robustness of nonlinear predictor feedback laws to time- and state-dependent delay perturbations", *Automatica*, to appear, 2013.

[6] Germani, A., C. Manes, and P. Pepe, "A New Approach to State Observation of Nonlinear Systems With Delayed Output", *IEEE Transactions on Automatic Control*, 47(1), 2002, 96-101.
[7] Karafyllis, I., "Stabilization By Means of Approximate Predictors for Systems with Delayed Input", *SIAM Journal on Control and Optimization*, 49(3), 2011, 1100-1123.
[8] Karafyllis, I., and Z.-P. Jiang, *Stability and Stabilization of Nonlinear Systems*, Springer-Verlag London (Series: Communications and Control Engineering), 2011.
[9] Karafyllis, I. and M. Krstic, "Nonlinear Stabilization under Sampled and Delayed Measurements, and with Inputs Subject to Delay and Zero-Order Hold", *IEEE Transactions on Automatic Control*, 57(5), 2012, 1141-1154.
[10] Karafyllis, I. and M. Krstic, "Predictor-Based Output Feedback for Nonlinear Delay Systems", arXiv:1108.4499v1 [math.OC].
[11] Karafyllis, I. and M. Krstic, "Numerical Schemes for Nonlinear Predictor Feedback", submitted to *Mathematics of Control, Signals and Systems* (see also arXiv:1211.1121 [math.OC]).
[12] Karafyllis, I. and M. Krstic, "On the Relation of Delay Equations to First-Order Hyperbolic Partial Differential Equations", submitted to *ESAIM Control, Optimisation and Calculus of Variations* (see also arXiv:1302.1128 [math.OC]).
[13] Khalil, H. K., *Nonlinear Systems*, $2^{nd}$ Edition, Prentice-Hall, 1996.
[14] Krstic, M., "Feedback Linearizability and Explicit Integrator Forwarding Controllers for Classes of Feedforward Systems", *IEEE Transactions on Automatic Control*, 49(10), 2004, 1668-1682.
[15] Krstic, M., "Lyapunov tools for predictor feedbacks for delay systems: Inverse optimality and robustness to delay mismatch", *Automatica*, 44(11), 2008, 2930-2935.
[16] Krstic, M., *Delay Compensation for Nonlinear, Adaptive, and PDE Systems*, Birkhäuser Boston, 2009.
[17] Krstic, M., "Input Delay Compensation for Forward Complete and Strict-Feedforward Nonlinear Systems", *IEEE Transactions on Automatic Control*, 55(2), 2010, 287-303.
[18] Mazenc, F., S. Mondie and R. Francisco, "Global Asymptotic Stabilization of Feedforward Systems with Delay at the Input", *IEEE Transactions on Automatic Control*, 49(5), 2004, 844-850.
[19] Mazenc, F., M. Malisoff and Z. Lin, "Further Results on Input-to-State Stability for Nonlinear Systems with Delayed Feedbacks", *Automatica*, 44(9), 2008, 2415-2421.
[20] Stuart, A.M. and A.R. Humphries, *Dynamical Systems and Numerical Analysis*, Cambridge University Press, 1998.
[21] Temam, R. *Infinite-Dimensional Dynamical Systems in Mechanics and Physics*, $2^{nd}$ Edition, Springer-Verlag, New York, 1997.